\documentclass[a4paper, 11pt, final]{article}

\usepackage{amsfonts,amsbsy}
\usepackage{enumerate}
\newtheorem{proposition}{Proposition}[section]
\newtheorem{theorem}[proposition]{Theorem}
\newtheorem{lemma}[proposition]{Lemma}

\newtheorem{definition}[proposition]{Definition}

\newtheorem{remark}[proposition]{Remark}

\newenvironment{proof}{\smallskip\noindent\emph{\textbf{Proof.}}\hspace{1pt}}%
{\hspace{-5pt}{\nobreak\quad\nobreak\hfill\nobreak$\square$\vspace{8pt}%
\par}\smallskip\goodbreak}

\newenvironment{proofof}[1]{\smallskip\noindent\emph{\textbf{Proof of #1.}}%
\hspace{1pt}}{\hspace{-5pt}{\nobreak\quad\nobreak\hfill\nobreak%
$\square$\vspace{8pt}\par}\smallskip\goodbreak}

\newcommand{\pint}[1]{\mathaccent23{#1}}

\newcommand{\C}[1]{\boldsymbol{\mathcal{C}^{#1}}}
\newcommand{\Cc}[1]{\boldsymbol{\mathcal{C}_c^{#1}}}
\newcommand{\modulo}[1]{{\left|#1\right|}}
\newcommand{\norma}[1]{{\left\|#1\right\|}}
\newcommand{\reali}{{\mathbb{R}}}
\newcommand{\rpic}{{\mathbb{R}_+}}
\newcommand{\rpis}{{\mathbb{R}_+^*}}
\newcommand{\naturali}{{\mathbb{N}}}
\newcommand{\BV}{\mathbf{BV}}

\renewcommand{\epsilon}{\varepsilon}
\renewcommand{\phi}{\varphi}
\renewcommand{\L}[1]{\mathbf{L^#1}}
\newcommand{\W}[2]{\mathbf{W^{#1,#2}}}

\newcommand{\Lloc}[1]{\mathbf{L^{#1}_{loc}}}
\newcommand{\tv}{\mathinner{\rm TV}}

\renewcommand{\div}{\mathinner{\rm div}}

\newcommand{\spt}{\mathop{\rm spt}}
\newcommand{\sgn}{\mathop{\rm sgn}}
\newcommand{\pt}{\partial}
\newcommand{\wto}{\rightharpoonup}
\newcommand{\Lip}{\mathrm{Lip}}

\renewcommand{\d}[1]{\mathinner{\mathrm{d}{#1}}}
\newcommand{\dx}{\d{x}}
\newcommand{\dt}{\d{t}}

\title{Control of the Continuity Equation \\ with a Non Local Flow}

\author{Rinaldo M. Colombo\footnote{Department of Mathematics, Brescia
  University, Via Branze 38, 25133 Brescia; Italy.}, Michael
Herty\footnote{RWTH Aachen University, Templergraben 55, 52056
  Aachen, Germany}, Magali Mercier\footnote{ Universit\'e de Lyon,
  Universit\'e Lyon1, \'Ecole Centrale de Lyon, INSA de Lyon, CNRS
  UMR 5208, Institut Camille Jordan; 43 blvd. du 11 Novembre 1918,
  F-69622 Villeurbanne Cedex, France}.}

\begin{document}

\maketitle

\begin{abstract}

  \noindent This paper focuses on the optimal control of weak (i.e.~in
  general \emph{non smooth}) solutions to the continuity equation with
  non local flow. Our driving examples are a supply chain model and an
  equation for the description of pedestrian flows.  To this aim, we
  prove the well posedness of a class of equations comprising these
  models. In particular, we prove the differentiability of solutions
  with respect to the initial datum and characterize its derivative. A
  necessary condition for the optimality of suitable integral
  functionals then follows.

  \medskip

\noindent\textit{2000~Mathematics Subject Classification:} 35L65,
49K20, 93C20

\medskip


\end{abstract}

\section{Introduction}
\label{sec:Intro}

We consider the continuity equation in $N$ space dimensions
\begin{equation}
  \label{eq:SCL}
  \left\{
    \begin{array}{l}
      \partial_t \rho + \div \left( \rho\; V (\rho) \right)
      = 0
      \\
      \rho(0,x) = \rho_o(x)
    \end{array}
  \right.
\end{equation}
with a \emph{non local} speed function $V$. This kind of equation
appears in numerous examples, a first one being the supply chain model
introduced in~\cite{ArmbrusterDegondRinghofer2006,
  ArmbrusterMarthalerRinghoferKempfJo2006}, where $V(\rho) = v \left(
  \int_0^1 \rho(x) \dx \right)$. Besides, this equation is very
similar to that obtained in a kinetic model of traffic,
see~\cite{BenzoniColomboG}. Another example comes from pedestrian
traffic, in which a reasonable model can be based on~(\ref{eq:SCL})
with the functional $V(\rho) = v(\rho\ast \eta) \,
\vec{v}(x)$. Throughout, our assumptions are modeled on these
examples.

The first question we address is that of the well posedness
of~(\ref{eq:SCL}). Indeed, we show in Theorem~\ref{thm:main}
that~(\ref{eq:SCL}) admits a unique local in time solution on a time
interval $I_{\mathrm{ex}}$. For all $t$ in $I_{\mathrm{ex}}$, we call
$S_t$ the nonlinear local semigroup that associates to the initial
condition $\rho_o$ the solution $S_t\rho_o$ of~(\ref{eq:SCL}) at time
$t$.  As in the standard case, $S_{t}$ turns out to be non expansive.

Then, we present a rigorous result on the G\^ ateaux differentiability
of the map $\rho_o \mapsto S_t {\rho_o}$, in any direction $r_o$ and for
all $t\in I_{\mathrm{ex}}$. Moreover, the G\^ ateaux derivative is
uniquely characterized as solution to the following linear Cauchy
problem, that can be obtained by linearising formally~(\ref{eq:SCL}):
\begin{equation}
  \label{eq:linear}
  \left\{
    \begin{array}{l}
      \pt_t r
      +
      \div \left( 
        r V(\rho) +  \rho \,  \mathrm{D}V(\rho)(r)
      \right)
      =
      0
      \\
      r(0,x)=r_o(x)\,.
    \end{array}
  \right.
\end{equation}
The well posedness of~(\ref{eq:linear}) is among the results of this
paper, see Proposition~\ref{prop:existence_r} below.

We stress here the difference with the well known standard
(i.e.~local) situation: the semigroup generated by a conservation law
is in general \emph{not} differentiable in $\L1$, not even in the
scalar 1D case, see~\cite[Section~1]{BressanGuerra}. To cope with
these issues, an entirely new differential structure was introduced
in~\cite{BressanGuerra}, and further developed in~\cite{Bianchini,
  BressanLewicka}, also addressing optimal control problems,
see~\cite{BressanShen2007, ColomboGroli2}. We refer
to~\cite{BouchutJames1998, BouchutJames1999,
  HertyGugatKlarLeugering2005b, U02, U03} for further results and
discussions about the scalar one--dimensional case. The presented
theories, however, seem not able to yield a \emph{``good''} optimality
criteria. On the one hand, several results deal only with smooth
solutions, whereas the rise of discontinuities is typical in
conservation laws. On the other hand, the mere definition of the shift
differential in the scalar 1D case takes alone about a page,
see~\cite[p.~89--90]{ColomboGroli2}.  Therefore, in the following we
postulate assumptions on the function $V$ which are satisfied in the
cases of the supply chain model and of the pedestrian model, but not
for general functions. To be more precise, we essentially require
below that $V$ is a \emph{non local} function,
see~(\ref{eq:TrulyNonLocal}).

Then, based on the differentiability results, we state a necessary
optimality condition.  We introduce a cost function $\mathcal{J}
\colon \C0 \left(I_{\mathrm{ex}}, \L1(\reali^N;\reali) \right) \to
\reali$ and, using the differentiability property given above, we find
a \emph{necessary condition} on the initial data $\rho_o$ in order to
minimize $\mathcal{J}$ along the solutions to~(\ref{eq:SCL})
associated to $\rho_o$.

We emphasize that all this is obtained within the framework of
\emph{non smooth} solutions, differently from many results in the
current literature that are devoted to
differentiability~\cite{Gugat2005}, control~\cite{Coron2007,
  CoronBastin2008} or optimal control~\cite{ColomboGuerraHertySachers}
for conservation laws, but limited to smooth solutions. Furthermore,
we stress that the present necessary conditions are obtained within
the functional setting typical of scalar conservation laws,
i.e.~within $\L1$ and $\L\infty$. No reflexivity property is ever
used.

The paper is organized as follows. In Section~\ref{sec:Main}, we state
the main results of this paper. The differentiability is proved in
Theorem~\ref{thm:Diff2} and applied to a control in supply chain
management in Theorem~\ref{thm:new}. The sections~\ref{sec:Supply}
and~\ref{sec:Ped} provide examples of models based on~(\ref{eq:SCL}),
and in Section~\ref{sec:Tech} we give the detailed proofs of our
results.

\section{Notation and Main Results}
\label{sec:Main}

\subsection{Existence of a Weak Solution to~(\ref{eq:SCL}) }

Denote $\reali_+ = \left[0, +\infty\right[$, $\rpis = \left]0,
  +\infty\right[$ and by $I$, respectively $I_*$ or $I_{\mathrm{ex}}$,
the interval $\left[0, T\right[$, respectively $\left[0, T_*\right[$
or $\left[ 0, T_{\mathrm{ex}}\right[$, for $T, T_*, T_{\mathrm{ex}} >
0$. Furthermore, we introduce the norms:
\begin{displaymath}
  \!\!\!\!\!\!
  \begin{array}{rcl@{\quad }rcl}
    \displaystyle
    \norma{v}_{\L\infty}
    & = &
    \displaystyle
    \sup_{x \in \reali^N} \norma{v(x)},
    &
    \displaystyle
    \norma{v}_{\W11}
    & = &
    \displaystyle
    \norma{v}_{\L1} + \norma{\nabla_xv}_{\L1}\,,
    \\
    \displaystyle
    \norma{v}_{\W2\infty}
    & = &
    \displaystyle
    \norma{v}_{\L\infty}+ \norma{\nabla_x v}_{\L\infty}
    + 
    \norma{\nabla_x^2 v}_{\L\infty},
    &
    \displaystyle
    \norma{v}_{\W1\infty}
    & = &
    \displaystyle
    \norma{v}_{\L\infty} + \norma{\nabla_xv}_{\L\infty}\,.
  \end{array}
  \!\!\!\!\!\!
\end{displaymath}
Let $V \colon \L1(\reali^N;\reali) \to \C2(\reali^N;\reali^N)$ be a
functional, not necessarily linear. A straightforward extension
of~\cite[Definition~1]{Kruzkov} yields the following definition of
weak solutions for~(\ref{eq:SCL}).
\begin{definition}
  \label{def:sol}
  A weak entropy solution to~(\ref{eq:SCL}) on $I_{\mathrm{ex}}$ is a
  bounded measurable map $\rho$ which is a Kru\v zkov solution to
  \begin{displaymath}
    \partial_t \rho + \div \left( \rho \, w(t,x) \right) = 0 \,,
    \qquad
    \mbox{ where} \qquad
    w(t,x) = \left(V \! \left(\rho(t) \right) \right) (x)  \,.
  \end{displaymath}
\end{definition}

\noindent In other words, for all $k \in \reali$ and for any test
function $\phi \in \Cc\infty (\pint{I}_{\mathrm{ex}} \times
\reali^N;\rpic)$
\begin{displaymath}
  \begin{array}{l}
    \displaystyle
    \int_0^{+\infty} \!\! \int_{\reali^N} \!\!
    \left[
      (\rho - k) \pt_t \phi
      +
      (\rho -k) \, {V \! \left(\rho(t)\right)(x)} \cdot \nabla_x \phi
      -
      \div \left( k \, {V \! \left(\rho(t)\right)(x)} \right) \phi
    \right]
    \\
    \displaystyle
    \qquad \qquad
    \times
    \sgn (\rho-k) \, \dx \, \dt
    \geq 0 
  \end{array}
\end{displaymath}
and there exists a set $\mathcal{E}$ of zero measure in $\rpic$ such
that for all $t \in I_{\mathrm{ex}}\setminus \mathcal{E}$ the function
$\rho$ is defined almost everywhere in $\reali^N$ and for any $\delta>
0$
\begin{displaymath}
  \lim_{t \to 0,\, t \in I_{\mathrm{ex}} \setminus \mathcal{E}}
  \int_{B(0,\delta)} \modulo{\rho(t,x) - \rho_o(x)} \mathrm{d}x
  =
  0 \,.
\end{displaymath}
The open ball in $\reali^N$ centered at $0$ with radius $\delta$ is
denoted by $B(0,\delta)$. Introduce the spaces
\begin{displaymath}
  \mathcal{X} 
  = 
  (\L1 \cap \L\infty \cap \BV) (\reali^N;\reali)
  \quad \mbox{ and } \quad 
  \mathcal{X}_{\alpha} 
  = 
  (\L1 \cap \BV) \bigl(\reali^N; [0,\alpha] \bigr)
  \mbox{ for } \alpha > 0
\end{displaymath}
both equipped with the $\L1$ distance. Obviously,
$\mathcal{X}_{\alpha} \subset \L\infty(\reali^N;\reali)$ for all
$\alpha > 0$.

We pose the following assumptions on $V$, all of which are satisfied
in the examples on supply chain and pedestrian flow as shown in
Section~\ref{sec:Supply} and Section~\ref{sec:Ped}, respectively.

\begin{description}
\item[(V1)] There exists a function $C\in \Lloc\infty(\rpic;\rpic)$
  such that for all $\rho\in \L1(\reali^N,\reali)$,
  \begin{eqnarray*}
    V(\rho) & \in & \L\infty (\reali^N; \reali^N)\,,
    \\
    \norma{\nabla_x V(\rho)}_{\L\infty(\reali^N; \reali^{N\times N})}
    & \leq &
    C(\norma{\rho}_{\L\infty(\reali^N;\reali)})\,,
    \\
    \norma{\nabla_x V(\rho)}_{\L1(\reali^N; \reali^{N\times N})} 
    & \leq & 
    C(\norma{\rho}_{\L\infty(\reali^N;\reali)}) \,,
    \\
    \norma{\nabla_x^2 V(\rho)}_{\L1(\reali^N; \reali^{N\times N\times N})} 
    & \leq & 
    C(\norma{\rho}_{\L\infty(\reali^N;\reali)})\,.
  \end{eqnarray*}
  There exists a function $C\in \Lloc\infty(\rpic;\rpic)$ such that
  for all $ \rho_1, \rho_2\in \L1(\reali^N,\reali)$
  \begin{eqnarray}
    \label{eq:TrulyNonLocal}
    \norma{V(\rho_1) - V(\rho_2)}_{\L\infty(\reali^N; \reali^N)}
    & \leq &
    C(\norma{\rho_1}_{\L\infty(\reali^N;\reali)}) \, 
    \norma{\rho_1 - \rho_2}_{\L1(\reali^N; \reali)}\,,
    \qquad
    \\
    \nonumber
    \norma{
      \nabla_x V(\rho_1) - \nabla_x V(\rho_2)
    }_{\L1(\reali^N; \reali^{N\times N})}
    & \leq &
    C(\norma{\rho_1}_{\L\infty(\reali^N;\reali)}) \, 
    \norma{\rho_1 - \rho_2}_{\L1(\reali^N; \reali)}
    \, .
  \end{eqnarray}
\item[(V2)] There exists a function $C\in \Lloc\infty(\rpic;\rpic)$
  such that for all $\rho\in \L1(\reali^N,\reali)$,
  \begin{displaymath}
    \norma{\nabla_x^2 V(\rho)}_{\L\infty(\reali^N;\reali^{N\times N\times N})}
    \leq
    C(\norma{\rho}_{\L\infty(\reali^N;\reali)})\, .
  \end{displaymath}
\item[(V3)] $V \colon \L1(\reali^N;\reali) \to \C3(\reali^N;\reali^N)$
  and there exists a function $C\in \Lloc\infty(\rpic;\rpic)$ such
  that for all $\rho\in \L1(\reali^N,\reali)$,
  \begin{displaymath}
    \norma{\nabla_x^3 V(\rho)}_{\L\infty(\reali^N;\reali^{N\times N\times N\times N})}
    \leq
    C(\norma{\rho}_{\L\infty(\reali^N;\reali)})\, .
  \end{displaymath}
\end{description}

\noindent Condition~(\ref{eq:TrulyNonLocal}) essentially requires that
$V$ be a \emph{non local} operator. Note that~\textbf{(V3)}
implies~\textbf{(V2)}. Existence of a solution to~(\ref{eq:SCL}) (at
least locally in time) can be proved under only
assumption~\textbf{(V1)}, see Theorem~\ref{thm:main}. The stronger
bounds on $V$ ensure additional regularity of the solution which is
required later to derive the differentiability properties, see
Proposition~\ref{prop:prop}.

\begin{theorem}
  \label{thm:main}
  Let~\textbf{(V1)} hold. Then, for all $\alpha,\beta>0$ with
  $\beta > \alpha$, there exists a time $ T(\alpha, \beta)>0$ such
  that for all $\rho_o \in \mathcal{X}_{\alpha}$,
  problem~(\ref{eq:SCL}) admits a unique solution $\rho \in \C0 \left(
    [0, T(\alpha,\beta)];\mathcal{X}_{\beta} \right)$ in the sense of
  Definition~\ref{def:sol}. Moreover,
  \begin{enumerate}
  \item $\norma{\rho(t)}_{\L\infty} \leq \beta$ for all $t \in [0,
    T(\alpha,\beta)]$.
  \item There exists a function $L \in \Lloc\infty(\rpic;\rpic)$ such
    that for all $\rho_{o,1}, \rho_{o,2}$ in $\mathcal{X}_{\alpha}$,
    the corresponding solutions satisfy, for all $t \in [0,
    T(\alpha,\beta)]$,
    \begin{displaymath}
      \norma{\rho_1(t) - \rho_2(t)}_{\L1}
      \leq
      L(t) \, \norma{\rho_{o,1} - \rho_{o,2}}_{\L1}
    \end{displaymath}
  \item There exists a constant $\mathcal{L} = \mathcal{L}(\beta)$
    such that for all $\rho_o \in \mathcal{X}_{\alpha}$, the
    corresponding solution satisfies for all $t \in [0,
    T(\alpha,\beta)]$
    \begin{displaymath}
      \tv\left(\rho(t) \right)
      \leq
      \left( 
        \tv(\rho_o) 
        + 
        \mathcal{L}   t \norma{\rho_o}_{\L\infty} 
      \right) e^{\mathcal{L}  t}
      \quad \mbox{ and } \quad
      \norma{\rho(t)}_{\L\infty}
      \leq
      \norma{\rho_o}_{\L\infty} \, e^{\mathcal{L}  t} \,.
    \end{displaymath}
  \end{enumerate}
\end{theorem}

\noindent The above result is \emph{local} in time. Indeed, as $t$
tends to $T(\alpha, \beta)$, the total variation of the solution may
well blow up. To ensure existence globally in time we need to
introduce additional conditions on $V$:
\begin{description}
\item[(A)] $V$ is such that for all $\rho \in \L1(\reali^N;\reali)$
  and all $x\in \reali^N$, $\left( \div V(\rho) \right) (x) \geq 0$.
\item[(B)] The function $C$ in~\textbf{(V1)} is bounded, i.e.~$C \in
  \L\infty(\rpic; \rpic)$.
\end{description}

\noindent Note that in the supply chain model discussed in
Section~\ref{sec:Supply}, condition~\textbf{(A)} applies.

\begin{lemma}
  \label{lem:Invariance}
  Assume all assumptions of Theorem~\ref{thm:main}. Let
  also~\textbf{(A)} hold. Then, for all $\alpha > 0$, the set
  $\mathcal{X}_\alpha$ is invariant for~(\ref{eq:SCL}), i.e.~if the
  initial datum $\rho_o$ satisfies
  $\norma{\rho_o}_{\L\infty(\reali^N;\reali)} \leq \alpha$, then,
  $\norma{\rho(t)}_{\L\infty(\reali^N;\reali)} \leq \alpha$ as long as
  the solution $\rho(t)$ exists.
\end{lemma}

Condition~\textbf{(B)}, although it does not guarantee the boundedness
of the solution, does ensure the global existence of the solution
to~(\ref{eq:SCL}).

\begin{theorem}
  \label{thm:V0}
  Let~\textbf{(V1)} hold. Assume moreover that~\textbf{(A)}
  or~\textbf{(B)} holds. Then, there exists a unique semigroup $S
  \colon \rpic \times \mathcal{X} \to \mathcal{X}$ with the following
  properties:
  \begin{enumerate}[(S1):]
  \item \label{it:s1} For all $\rho_o \in \mathcal{X}$, the orbit $t
    \mapsto S_t \rho_o$ is a weak entropy solution to~(\ref{eq:SCL}).
  \item \label{it:s4}$S$ is $\L1$-continuous in time, i.e.~for all
    $\rho_o \in \mathcal{X}$, the map $t \mapsto S_t \rho_o$ is in
    $\C0 (\reali_+; \mathcal{X})$.
  \item \label{it:s3}$S$ is $\L1$-Lipschitz with respect to the
    initial datum, i.e.~for a suitable positive $L \in \Lloc\infty
    (\rpic;\rpic)$, for all $t \in \reali_+$ and all $\rho_1,\rho_2
    \in \mathcal{X}$,
    \begin{displaymath}
      \norma{S_t \rho_1 - S_t \rho_2}_{\L1(\reali^N;\reali)}
      \leq
      L(t) \, \norma{\rho_1 - \rho_2}_{\L1(\reali^N;\reali)} \,.
    \end{displaymath}
  \item \label{it:s2} There exists a positive constant $\mathcal{L}$
    such that for all $\rho_o \in \mathcal{X}$ and all $t \in \rpic$,
    \begin{displaymath}
      \tv\left(\rho(t) \right)
      \leq 
      \left( 
        \tv(\rho_o) 
        + 
        \mathcal{L}  t \, \norma{\rho_o}_{\L\infty(\reali^N;\reali)} 
      \right) e^{\mathcal{L}t} \,.
    \end{displaymath}
  \end{enumerate}
\end{theorem}

\noindent Higher regularity of the solutions of~(\ref{eq:SCL}) can be
proved under stronger bounds on $V$.

\begin{proposition}
  \label{prop:prop}
  Let~\textbf{(V1)} and~\textbf{(V2)} hold. With the same
  notations as in Theorem~\ref{thm:main}, if $\rho_o \in
  \mathcal{X}_\alpha$, then
  \begin{displaymath}
    \begin{array}{rclcrcl@{\quad}rcl}
      \rho_o & \in & (\W11 \cap \L\infty)(\reali^N;\reali)
      & \Longrightarrow &
      \forall t & \in & [0, T(\alpha,\beta)],
      & \rho(t) & \in & \W11(\reali^N;\reali)\, ,
      \\ 
      \rho_o & \in & \W1\infty (\reali^N;\reali)
      & \Longrightarrow &
      \forall t & \in & [0, T(\alpha,\beta)], 
      & \rho(t) & \in & \W1\infty(\reali^N;\reali)\, ,
    \end{array}
  \end{displaymath}
  and there exists a positive constant $C = C(\beta)$ such that
  \begin{displaymath}
    \norma{\rho(t)}_{\W11} 
    \leq
    e^{2Ct} \, \norma{\rho_o}_{\W11}
    \quad \mbox{ and } \quad
    \norma{\rho(t)}_{\W1\infty}
    \leq
    e^{2Ct} \, \norma{\rho_o}_{\W1\infty} \,.
  \end{displaymath}
  Furthermore, if $V$ also satisfies~\textbf{(V3)}, then
  \begin{eqnarray*}
    \rho_o \in (\W21 \cap \L\infty)(\reali^N;[\alpha,\beta])
    &\Longrightarrow &
    \forall t\in  [0, T(\alpha,\beta)], \quad
    \rho(t)\in \W21(\reali^N;\reali)\,
  \end{eqnarray*}
  and for a suitable non--negative constant $C = C(\beta)$, we have
  the estimate
  \begin{eqnarray*}
    \norma{\rho(t)}_{\W21} 
    & \leq &
    e^{Ct}(2e^{Ct}-1)^2 \, \norma{\rho_o}_{\W21} \,.
  \end{eqnarray*}
\end{proposition}

The proofs are deferred to Section~\ref{sec:Tech}.

\subsection{Differentiability}
\label{sec:Diff}

This section is devoted to the differentiability of the semigroup $S$
(defined in Theorem~\ref{thm:main}) with respect to the initial datum
$\rho_{o}$, according to the following notion. Recall first the
following definition.

\begin{definition}
  \label{def:weak diff}
  A map $F \colon \L1(\reali^N; \reali) \to \L1(\reali^N; \reali)$ is
  \emph{strongly $\L1$ G\^ateaux differentiable in any direction} at
  $\rho_o \in \L1(\reali^N; \reali)$ if there exists a continuous
  linear map $DF(\rho_o) \colon \L1(\reali^N;\reali) \to
  \L1(\reali^N;\reali)$ such that for all $r_o \in
  \L1(\reali^N;\reali)$ and for any real sequence $(h_n)$ with $h_n
  \to 0$,
  \begin{displaymath}
%     \frac{F(\rho_o + h_n r_o) - F(\rho_o)}{h_n}
%     & \stackrel{n \to \infty}{\wto} &
%     DF(\rho_o) (r_o)
%     \quad \hbox{ weakly in } \L1 
%     \mbox{, respectively}
%     \\
    \frac{F(\rho_o + h_n r_o) - F(\rho_o)}{h_n}
    \stackrel{n\to \infty}{\to} 
    DF(\rho_o) (r_o)
    \quad \hbox{ strongly in } \L1 
    .
  \end{displaymath}
\end{definition}

Besides proving the differentiability of the semigroup, we also
characterize the differential. Formally, a sort of first order
expansion of~(\ref{eq:SCL}) with respect to the initial datum can be
obtained through a standard linearization procedure, which
yields~(\ref{eq:linear}). Now, we rigorously show that the derivative
of the semigroup in the direction $r_o$ is indeed the solution
to~(\ref{eq:linear}) with initial condition $r_o$. To this aim, we
need a forth and final condition on $V$.

\begin{description}
\item[(V4)] $V$ is Fr\'echet differentiable as a map $\L1 (\reali^N;
  \rpic) \to \C2 (\reali^N; \reali^N)$ and there exists a function $K
  \in \Lloc\infty(\rpic; \rpic)$ such that for all $\rho \in \L1 (
  \reali^N; \rpic)$, for all $r\in \L1(\reali^N;\reali)$,
  \begin{eqnarray*}
    \norma{V(\rho + r) - V(\rho) -
      DV(\rho) (r)}_{\W2\infty}
    &\leq &
    K
    \left(
      \norma{\rho}_{\L\infty}
      +
      \norma{\rho+r}_{\L\infty}
    \right)
    \, \norma{r}_{\L1}^2 \,,
    \\
    \norma{DV(\rho)(r)}_{\W2\infty} 
    &\leq & 
    K \left( \norma{\rho}_{\L\infty} \right) \,
    \norma{r}_{\L1}\,.
  \end{eqnarray*}
\end{description}

\noindent Consider now system~(\ref{eq:linear}), where $\rho \in \C0
(I_{\mathrm{ex}}, \mathcal{X})$ is a given function. We introduce a
notion of solution for~(\ref{eq:linear}) and give conditions which
guarantee the existence of a solution.

\begin{definition}
  Fix $r_o \in \L\infty(\reali^N;\reali)$. A function $r \in \L\infty
  \bigl( I_{\mathrm{ex}}; \Lloc1(\reali^N; \rpic) \bigr)$ bounded,
  mesurable and right continuous in time, is a \emph{weak solution}
  to~(\ref{eq:linear}) if for any test function $\phi \in \Cc\infty
  (\pint{I}_{\mathrm{ex}} \times \reali^N;\reali)$
  \begin{equation}
    \label{eq:defWeak}
    \begin{array}{l}
      \displaystyle
      \int_0^{+\infty} \int_{\reali^N}
      \left[
        r \, \partial_t \phi
        +
        r \, a(t,x)\cdot \nabla_x \phi - \div b(t,x) \phi
      \right]
      \, \dx \, \dt
      = 0\, , \quad      \mbox{ and}
      \\[10pt]
      r(0) = r_o\,  \quad \mbox{ a.e.~in } \; \reali^N,
    \end{array}
  \end{equation}
  where $a= V(\rho)$ and $b=\rho DV(\rho)(r)$.
\end{definition}

\noindent We now extend the classical notion of Kru\v zkov solution to
the present non local setting.

\begin{definition}
  Fix $r_o \in \L\infty(\reali^N; \rpic)$. A function $r \in \L\infty
  \bigl( I_{\mathrm{ex}}; \Lloc1(\reali^N; \rpic) \bigr)$ bounded,
  mesurable and right continuous in time, is a \emph{Kru\v zkov
    solution to the nonlocal problem~(\ref{eq:linear})} if it is a
  Kru\v{z}kov solution to
  \begin{equation}
    \label{eq:ab}
    \left\{
      \begin{array}{l}
        \partial_t r +\div \left( r \, a(t,x) + b(t,x) \right) =0
        \\
        r(0,x) = r_o(x)
      \end{array}
    \right.
  \end{equation}
  where $a= V(\rho)$ and $b=\rho DV(\rho)(r)$.
\end{definition}
In other words, $r$ is a Kru\v{z}kov solution to~(\ref{eq:linear}) if
for all $k \in \reali$ and for any test function $\phi \in \Cc\infty
(\pint{I}_{\mathrm{ex}} \times \reali^N;\rpic)$
\begin{displaymath}
  \!\!\!\!\!\!\!
  \begin{array}{l}
    \displaystyle
    \int_0^{+\infty} \!\!\! \int_{\reali^N} \!\!
    \left[
      (r - k) \pt_t \phi
      +
      (r -k) V(\rho) \cdot \nabla_x \phi
      -
      \div \left( k V(\rho) + \rho DV(\rho)r \right) \phi
    \right]
    \sgn (r-k) \dx \dt
    \geq 0
    \\[5pt]
    \mbox{ and}  
    \quad 
    \displaystyle
    \lim_{t\to 0^+} \int_{B(0,\delta)} 
    \modulo{r(t) - r_o} \dx = 0 
    \quad \mbox{ for all } \delta > 0 \,.
  \end{array}
  \!\!\!\!
\end{displaymath}
Condition~\textbf{(V4)} ensures that if $\rho \in
\W11(\reali^N;\reali)$, then $DV(\rho)(r)\in \C2(\reali^N; \reali^N)$
and hence for all $t \geq 0$, the map $x \mapsto \rho(t,x) \,
DV\left(\rho(t)\right) r(t,x)$ is in $\W11(\reali^N;\reali)$, so that
the integral above is meaningful.

\begin{proposition}
  \label{prop:existence_r}
  Let~\textbf{(V1)} and~\textbf{(V4)} hold. Fix $\rho \in \C0
  \bigl(I_{\mathrm{ex}}; (\W1\infty\cap \W11)(\reali^N; \reali)
  \bigr)$. Then, for all $r_o \in (\L1\cap\L\infty)(\reali^N;\reali)$
  there exists a unique weak entropy solution to~(\ref{eq:linear}) in
  $\L\infty \bigl( I_{\mathrm{ex}}; \L1(\reali^N;\reali) \bigr)$
  continuous from the right in time, and for all time $t \in
  I_{\mathrm{ex}}$, with $C =
  C\left(\norma{\rho}_{\L\infty([0,t]\times\reali^N;\reali)}\right)$
  as in~\textbf{(V1)} and $K =
  K\left(\norma{\rho}_{\L\infty([0,t]\times\reali^N;\reali)}\right)$
  as in~\textbf{(V4)}
  \begin{eqnarray*}
    \norma{r(t)}_{\L1}
    & \leq &
    e^{Kt\norma{\rho}_{\L\infty(I;\W11)}} \, e^{Ct} \, 
    \norma{r_o}_{\L1}
    \\
    \norma{r(t)}_{\L\infty}
    & \leq &
    e^{Ct} \, \norma{r_o}_{\L\infty} 
    + 
    K \, t \, e^{2Ct} \, e^{Kt\norma{\rho}_{\L\infty(I;\W11)}} \,
    \norma{\rho}_{\L\infty(I;\W1\infty)} \, \norma{r_o}_{\L1} \,.
  \end{eqnarray*}
  If~\textbf{(V2)} holds, $\rho \in \L\infty \bigl(I_{\mathrm{ex}};
  (\W1\infty \cap \W21) (\reali^N;\reali) \bigr)$ and $r_o \in (\W11
  \cap \L\infty)(\reali^N;\reali)$, then for all $t \in I_{\mathrm{ex}}
  $, $r(t) \in \W11(\reali^N;\reali)$ and
  \begin{displaymath}
    \norma{r(t)}_{\W11}
    \leq 
    (1+C't) \, e^{2C't} \, \norma{r_o}_{\W11} 
    + 
    K t (1 +  Ct) \, e^{4C't} \,
    \norma{r_o}_{\L1} \,
    \norma{\rho}_{\L\infty(I;\W21)}
    \,.
  \end{displaymath}
  where $C' = \max \bigl\{ C, K \norma{\rho}_{\L\infty(I_{\mathrm{ex}};
      \W21(\reali^N;\reali))} \bigr\}$. Furthermore, full continuity
  in time holds: $r \in \C0(I_{\mathrm{ex}};\L1(\reali^N;\reali))$.
\end{proposition}

With these tools, we can now state a theorem about the weak G\^ateaux
differentiability.

\begin{theorem}
  \label{thm:Diff}
  Let~\textbf{(V1)} and~\textbf{(V4)} hold. Let $\rho_o \in (\W1\infty
  \cap \W11)(\reali^N; \reali)$, and denote $T_{\mathrm{ex}}$ the time
  of existence for the solution of~(\ref{eq:SCL}). Then, for all time
  $t\in I_{ex}$, for all $r_o\in \mathcal{X}$ and all sequences
  $(h_n)_{n \in \naturali}$ converging to 0, there exists a
  subsequence of $\left(\frac{1}{h_n} \left( S_t( \rho_o + h_n
      r_o)-S_t(\rho_o) \right) \right)_{n\in \naturali}$ that
  converges weakly in $\L1$ to a weak solution of~(\ref{eq:linear}).
\end{theorem}

This theorem does not guarantee the uniqueness of this kind of \textit{weak $\L1$ G\^ateaux
derivative}.  Therefore, we consider the following stronger hypothesis,
under which we derive a result of strong G\^ateaux differentiability
and uniqueness of the derivative.
\begin{description}
\item[(V5)] There exists a function $K \in \Lloc\infty(\rpic;\rpic)$
  such that $\forall \rho, \tilde \rho \in \L1(\reali^N;\reali)$
  \begin{displaymath}
    \norma{
      \div \left( V(\tilde\rho) -  V(\rho) -
        DV(\rho) (\tilde \rho-\rho)\right)}_{\L1}
    \leq
    K
    \left(
      \norma{\rho}_{\L\infty}
      +
      \norma{\tilde\rho}_{\L\infty}
    \right) 
    \left( \norma{\tilde \rho - \rho}_{\L1} \right)^2 
  \end{displaymath}
  and the map $r \to \div DV(\rho)(r)$ is a bounded linear operator on
  $\L1(\reali^N; \reali) \to \L1(\reali^N;\reali)$, i.e.~$\forall
  \rho, r \in \L1(\reali^N;\reali)$
  \begin{displaymath}
    \norma{\div \left(DV(\rho)(r) \right)}_{\L1(\reali^N;\reali)}
    \leq
    K \left( \norma{\rho}_{\L\infty(\reali^N;\reali)} \right) \,
    \norma{r}_{\L1(\reali^N;\reali)}\,.
  \end{displaymath}
\end{description}

\begin{theorem}
  \label{thm:Diff2}
  Let~\textbf{(V1)}, \textbf{(V3)}, \textbf{(V4)} and~\textbf{(V5)}
  hold. Let $\rho_o \in (\W1\infty \cap \W21 ) (\reali^N; \reali)$,
  $r_o \in (\W11 \cap \L\infty) (\reali^N; \reali)$, and denote
  $T_{\mathrm{ex}}$ the time of existence of the solution
  of~(\ref{eq:SCL}) with initial condition $\rho_o$. Then, for all
  time $t\in I_{\mathrm{ex}}$ the local semigroup defined in
  Theorem~\ref{thm:main} is strongly $\L1$ G\^{a}teaux differentiable
  in the direction $r_o$. The derivative $DS_t(\rho_o)(r_o)$ of
  $S_t$ at $\rho_o$ in the direction $r_o$ is 
  \begin{displaymath}
    DS_t(\rho_o)(r_o) = \Sigma_t^{\rho_o}(r_o) \, .
  \end{displaymath}
  where $\Sigma^{\rho_o}$ is the linear application generated by the
  Kru\v{z}kov solution to~(\ref{eq:linear}), where $\rho = S_t
  \rho_o$, then for all $t\in I_{\mathrm{ex}}$. 
%Moreover, $DS_t(\rho_o) \colon \L1 \to\L1$ is linear and continuous.
\end{theorem}

\subsection{Necessary Optimality Conditions for Problems Governed
  by~(\ref{eq:SCL})}
\label{subsec:Optimal}

Aiming at necessary optimality conditions for non linear functionals
defined on the solutions to~(\ref{eq:SCL}), we prove the following
chain rule formula.

\begin{proposition}
  \label{prop:case} Let $T>0$ and $I=[0,T[$.  Assume that $f\in
  \C{1,1}(\reali;\rpic)$, $\psi \in \L\infty(I \times \reali^N;
  \reali)$ and that $S \colon I \times (\L1 \cap \L\infty)
  (\reali^N;\reali) \to (\L1 \cap \L\infty) (\reali^N;\reali)$ is
  strongly $\L1$ G\^ateaux differentiable. For all $t \in I$, let
  \begin{equation}
    \label{eq:J}
    J (\rho_{o} ) 
    = 
    \int_{\reali^N} 
    f \left( S_t \rho_o \right) \, \psi(t,x)
    \dx \,.
  \end{equation}
  Then, $J$ is strongly $\L\infty$ G\^ateaux differentiable in any
  direction $r_o \in (\W11 \cap \L\infty) (\reali^N;\reali)$. Morever,
  \begin{displaymath}
    DJ(\rho_o)(r_o) 
    =
    \int_{\reali^N} 
    f'(S_t\rho_o) \, \Sigma_t^{\rho_o} (r_o) \psi(t,x) \, \dx \, . 
  \end{displaymath}
\end{proposition}

\begin{proof}
  Since $\modulo{ f(\rho_h) - f(\rho) - f'(\rho) (\rho_h-\rho)} \leq
  \Lip(f') \, \modulo{\rho_h-\rho}^2$, we have
  \begin{eqnarray*}
    & &
    \modulo{
      \frac{J(\rho_o+hr_o)-J(\rho_o)}{h}
      -
      \int_{\reali^N}
      f'(S_t\rho_o) \, DS_t(\rho_o)(r_o) \, \psi(t,x) \, \dx}
    \\
    & \leq &
    \int_{\reali^N}  \modulo{f'(S_t\rho_o)}
    \modulo{
      \frac{S_t(\rho_o+hr_o)-S_t(\rho_o)}{h}
      -
      DS_t(\rho_o)(r_o)
    }
    \, \modulo{\psi(t,x)} \, \dx
    \\
    & &
    + 
    \Lip(f') \, \frac{1}{\modulo{h}} \int_{\reali^N} 
    \modulo{S_t(\rho_o + h r_o) - S_t(\rho_o)}^2 \, 
    \modulo{\psi(t,x)}
    \, \dx \,.
  \end{eqnarray*}
  The strong G\^ateaux differentiability of $S_t$ in $\L1$ then yields
  \begin{displaymath}
    \int_{\reali^N}
    \modulo{f'(S_t\rho_o)} \,
    \modulo{
      \frac{S_t(\rho_o + h r_o) - S_t(\rho_o)}{h}
      -
      DS_t(\rho_o)(r_o) 
    }
    \, \modulo{\psi(t,x)} \, \dx  = o(1)
    \quad \mbox{ as } h \to 0
  \end{displaymath}
  thanks to $S_t\rho_o \in \L\infty$ and to the local boundedness of
  $f'$.  Furthermore,
  \begin{eqnarray*}
    S_t(\rho_o) , S_t(\rho_o+hr_o) & \in & \L\infty
    \\
    \frac{1}{h} \left( S_t(\rho_o + h r_o) - S_t (\rho_o) \right) 
    & \stackrel{h\to 0}{\longrightarrow} &
    D S_t (\rho_o) (r_o) \mbox{ pointwise a.e.}
    \\
    S_t(\rho_o + h r_o) - S_t (\rho_o) 
    & \stackrel{h\to 0}{\longrightarrow} &
    0 \mbox{ pointwise a.e.}
  \end{eqnarray*}
  the Dominated Convergence Theorem ensures that the higher order term
  in the latter expansion tend to $0$ as $h \to 0$.
\end{proof}

The above result can be easily extended. First, to more general (non
linear) functionals $J(\rho_o) = \mathcal{J} (S_t \rho_o)$, with
$\mathcal{J}$ satisfying
\begin{description}
\item[(J)] $\mathcal{J} \colon \mathcal{X} \to \rpic$ is Fr\'echet
  differentiable: for all $\rho \in \mathcal{X}$ there exists a
  continuous linear application $D\mathcal{J}(\rho) \colon \mathcal{X}
  \to \reali$ such that for all $\rho, r \in \mathcal{X}$:
  \begin{displaymath}
    \modulo{
      \frac{\mathcal{J}(\rho+hr) - \mathcal{J}(\rho)}{h}
      -
      D\mathcal{J}(\rho)(r)
    } 
    \stackrel{h\to 0}{\longrightarrow} 0\, .
  \end{displaymath}
\end{description}
\noindent Secondly, to functionals of the type
\begin{displaymath}
  J(\rho_o) 
  =
  \int_0^T \int_{\reali^N}
  f(S_t \rho_o) \, \psi(t,x) \, \dx \, \dt
  \quad \mbox{ or } \quad
  J(\rho_o) 
  = 
  \int_0^T \mathcal{J} (S_t \rho_o) \, \dt \,.
\end{displaymath}
This generalization, however, is immediate and we omit the details.

Once the differentiability result above is available, a necessary
condition of optimality is straightforward.

\begin{proposition}
  \label{prop:optimal}
  Let $f\in \C{1,1}(\reali;\rpic)$ and $\psi \in
  \L\infty(I_{\mathrm{ex}} \times \reali^N; \reali)$. Assume that $S
  \colon I \times (\L1 \cap \L\infty) (\reali^N;\reali) \to (\L1 \cap
  \L\infty) (\reali^N;\reali)$ is strongly $\L1$ G\^ateaux
  differentiable. Define $J$ as in~(\ref{eq:J}). If $\rho_o \in (\L1
  \cap \L\infty) (\reali^N;\reali)$ solves the problem
  \begin{displaymath}
    % \label{eq:optimal}
    \textrm{find } \min_{\rho_{o}} \mathcal{J}( \rho  ) 
    \mbox{ subject to } \{\rho 
    \textrm{ is solution to~(\ref{eq:SCL})}\}.
  \end{displaymath}
  then, for all $r_o \in (\L1 \cap \L\infty) (\reali^N; \reali)$
  \begin{equation}
    \label{eq:last}
    \int_{\reali^N} f'(S_t \rho_o) \, \Sigma_t^{\rho_o} r_o \, \psi(t,x)
    \, \dx 
    =
    0 \,.
  \end{equation}
\end{proposition}

\section{Demand Tracking Problems for Supply Chains}
\label{sec:Supply}

Recently, Armbruster et
al.~\cite{ArmbrusterMarthalerRinghoferKempfJo2006}, introduced a
continuum model to simulate the average behavior of highly re-entrant
production systems at an aggregate level appearing, for instance, in
large volume semiconductor production line. The factory is described
by the density of products $\rho(x,t)$ at stage $x$ of the production
at a time $t$.  Typically, see~\cite{Ag76,
  ArmbrusterMarthalerRinghoferKempfJo2006, Karma89}, the production
velocity $V$ is a given smooth function of the total load $\int_0^1
\rho(t,x) \dx$, for example
\begin{equation}
  \label{v_eq}
  v(u) = v_{\max} / (1+u)
  \quad \mbox{ and } \quad
  V(\rho) = v \left( \int_0^1 \rho(t,s) \, ds \right) \,.
\end{equation}
The full model, given by~(\ref{eq:SCL})--(\ref{v_eq}) with $N=1$, fits
in the present framework.

\begin{proposition}
  \label{prop:OptimalSupply}
  Let $v \in \C1\left( [0,1]; \reali\right)$. Then, the functional $V$
  defined as in~(\ref{v_eq}) satisfies~\textbf{(A)}, \textbf{(V1)},
  \textbf{(V2)}, \textbf{(V3)}.  Moreover, if $v \in \C2\left( [0,1];
    \reali \right)$, then $V$ satisfies also~\textbf{(V4)}
  and~\textbf{(V5)}.
\end{proposition}

\noindent The proof is deferred to Paragraph~\ref{subsec:ApplProofs}.

The supply chain model with $V$ given by~(\ref{v_eq})
satisfies~\textbf{(V1)} to~\textbf{(V5)} and~\textbf{(A)}. Therefore,
Theorem~\ref{thm:V0} applies and, in particular, the set $[0,1]$ is
invariant yielding global well posedness. By Theorem~\ref{thm:Diff2},
the semigroup $S_{t} \rho_{o}$ is G\^{a}teaux differentiable in any
direction $r_{o}$ and the differential is given by the solution
to~(\ref{eq:linear}).

Note that the velocity is constant across the entire system at any
time. In fact, in a real world factory, all parts move through the
factory with the same speed. While in a serial production line, speed
through the factory is dependent on all items and machines downstream,
in a highly re-entrant factory this is not the case. Since items must
visit machines more than once, including machines at the beginning of
the production process, their speed through factory is determined by
the total number of parts both upstream and downstream from them. Such
re-entrant production is characteristic for semiconductor fabs.
Typically, the output of the whole factory over a longer timescale,
e.g.~following a seasonal demand pattern or ramping up or down a new
product, can be controlled by prescribing the inflow density to a
factory $\rho(t,x=0) = \lambda(t)$. The influx should be chosen in
order to achieve either of the following objective
goals~\cite{ArmbrusterMarthalerRinghoferKempfJo2006}:
\begin{enumerate}
\item[(1)] Minimize the mismatch between the outflux and a demand rate
  target $d(t)$ over a fixed time period (demand tracking problem).
  This is modelled by the cost functional $\frac{1}{2} \int_{0}^T
  \left( d(t) - \rho(1,t) \right)^2 \, dt. $
\item[(2)] Minimize the mismatch between the total number of parts
  that have left the factory and the desired total number of parts
  over a fixed time period $d(t)$. The backlog of a production system
  at a given time $t$ is defined as the total number of items that
  have been demanded minus the total number of items that have left
  the factory up to that time. Backlog can be negative or positive,
  with a negative backlog corresponding to overproduction and a
  positive backlog corresponding to a shortage.  This problem is
  modeled by $ \frac{1}{2} \int_0^T \left( \int_{0}^t d(\tau) u-
    \rho(1,\tau) d\tau \right)^2 dt$.
\end{enumerate}
In both cases we are interested in the influx $\lambda(t)$. A
numerical integration of this problem has been studied
in~\cite{LaMarcaHertyArmbrusterHertyRinghofer2008}. In order to apply
the previous calculus we reformulate the optimization problem for the
influx density $\lambda(t) = \rho_{o}(-t)$ where $\rho_{0}$ is the
solution to a minimization problem for
\begin{equation}
  \label{eq:Jsupply}
  \begin{array}{rcl}
    J_{1}( \rho_{o}) 
    & = &
    \displaystyle 
    \frac{1}2 \int_0^1 ( d(x) -  S_{T}\rho_{o}(x) )^2 \d{x}
    \\
    J_{2}( \rho_{o}) 
    & = &
    \displaystyle 
    \frac{1}2 
    \int_0^1 
    \left( 
      \int_{0}^x \left( d(\xi) -  S_{T}\rho_{o}(\xi)
      \right) \d\xi 
    \right)^2 
    \d{x}\,,
  \end{array}
\end{equation}
respectively, where $S_{t}\rho_{o}$ is the solution to~(\ref{eq:SCL})
and~(\ref{v_eq}). Clearly, $J_{1}$ and $J_{2}$ satisfy the assumptions
imposed in the previous section. The assertions of
Proposition~\ref{prop:optimal} then state necessary optimality
conditions, which we summarize in the theorem below.

\begin{theorem}
  \label{thm:new}
  Let $T > 0$ be given. Let the assumptions of
  Proposition~\ref{prop:OptimalSupply} hold. Let $\rho_{o} \in
  (\W1\infty\cap \W21)(\reali; \reali)$ be a minimizer of $J_1$ as
  defined in~(\ref{eq:Jsupply}), with $S$ being the semigroup
  generated by~(\ref{eq:SCL})--(\ref{v_eq}).  Then, for all $r_{o} \in
  (\W1\infty\cap \W21)(\reali;\reali)$ we have
  \begin{eqnarray*}
    \int_{0}^1 \left( d(x) - \rho(T,x) \right) r(T,x) \dx 
    & = &
    0\,, \quad \mbox{where}
    \\
    \partial_{t} r 
    + 
    \partial_{x} 
    \left( 
      v_{\max} 
      \frac{\left( 
          r \int_{0}^1 \rho \dx + \rho \int_{0}^1 r \dx 
        \right)}{\left( 
          1 + \int_{0}^1 \rho \dx \right)^2} 
    \right)
    & = &
    0\,, \quad 
    r(0,x) = r_{o}(x) \,.
  \end{eqnarray*}  
\end{theorem}

\noindent The latter Cauchy problem is in the form~(\ref{eq:linear})
and Proposition~\ref{prop:existence_r} proves its well posedness. The
latter proof is deferred to Paragraph~\ref{subsec:ApplProofs}.

\section{A Model for Pedestrian Flow}
\label{sec:Ped}

Macroscopic models for pedestrian movements are based on the
continuity equation, see~\cite{ColomboRosini1, Hughes2002,
  PiccoliTosin}, possibly together with a second equation, as
in~\cite{CosciaCanavesio}. In these models, pedestrians are assumed to
instantaneously adjust their (vector) speed according to the crowd
density at their position. The analytical construction in
Section~\ref{sec:Main} allows to consider the more realistic situation
of pedestrian deciding their speed according to the local mean density
at their position. We are thus led to consider~(\ref{eq:SCL}) with
\begin{equation}
  \label{eq:Ped}
  V(\rho) =  v( \rho * \eta) \, \vec{v}
\end{equation}
where
\begin{equation}
  \label{eq:eta}
  \eta \in \Cc2 \bigl( \reali^2; [0,1] \bigr)
  \mbox{ has support } \spt \eta \subseteq B(0,1)
  \mbox{ and }\norma{\eta}_{\L1} = 1 \,,
\end{equation}
so that $(\rho * \eta) (x)$ is an average of the values attained by
$\rho$ in $B(x,1)$. Here, $\vec v = \vec v(x)$ is the given direction
of the motion of the pedestrian at $x \in \reali^2$. Then, the
presence of boundaries, obstacles or other geometric constraint can be
described through $\vec v$, see~\cite{ColomboFacchiMaterniniRosini,
  PiccoliTosin}.

Note here that condition~\textbf{(A)} is unphysical, for it does not
allow any increase in the crowd density. Hence, for this example we
have only a local in time solution by Theorem~\ref{thm:main}.

As in the preceding example, first we state the hypotheses that
guarantee assumptions~\textbf{(V1)} to~\textbf{(V5)}.

\begin{proposition}
  \label{prop:ped}
  Let $V$ be defined in~(\ref{eq:Ped}) and $\eta$ be as
  in~(\ref{eq:eta}).
  \begin{enumerate}
  \item If $v \in \C2\left( \reali; \reali \right)$ and $\vec{v} \in
    (\C2\cap \W21) (\reali^2; \mathbb{S}^1)$, then $V$
    satisfies~\textbf{(V1)} and~\textbf{(V2)}.
  \item If moreover $v \in \C3(\reali;\reali)$, $\vec{v} \in
    \C3(\reali^2; \reali^2)$ and $\eta \in \C3(\reali^2;\reali)$ then
    $V$ satisfies~\textbf{(V3)}.
  \item If moreover $v\in \C4(\reali;\reali)$, $\vec{v}\in
    \C2(\reali^2;\reali^2)$ and $\eta\in \C2(\reali^2;\reali)$, then
    $V$ satisfies~\textbf{(V4)} and~\textbf{(V5)}.
  \end{enumerate}
\end{proposition}

\noindent The proof is deferred to Paragraph~\ref{subsec:ApplProofs}.

A typical problem in the management of pedestrian flows consists in
keeping the crowd density $\rho(t,x)$ below a given threshold, say
$\hat\rho$, in particular in a sensible compact region $\Omega$. To
this aim, it is natural to introduce a cost functional of the type
\begin{equation}
  \label{eq:Jped}
  J (\rho_o)
  =
  \int_0^T \int_{\reali^N}
  f\left( S_t \rho_o(x) \right) \, \psi(t,x) \dx \dt
\end{equation}
where
\begin{description}
\item[(f)] $f \in \C{1,1}(\reali;\rpic)$, $f(\rho)=0$ for $\rho \in
  [0,\hat\rho]$, $f(\rho) > 0$ and $f'(\rho) > 0$ for $\rho >
  \hat\rho$.
\item[($\boldsymbol{\psi}$)] $g \in \C\infty (\reali^N;[0,1])$, with
  $\spt g = \Omega$, is a smooth approximation of the characteristic
  function of the compact set $\Omega$, with $\pint{\Omega} \neq
  \emptyset$.
\end{description}

\noindent Paragraph~\ref{subsec:Optimal} then applies, yielding the
following necessary condition for optimality.

\begin{theorem}
  \label{thm:newped}
  Let $T > 0$ and the assumptions of 1.--3.~in
  Proposition~\ref{prop:ped} hold, together with~\textbf{(f)}
  and~\textbf{($\boldsymbol{\psi}$)}. Let $\rho_{o} \in (\W1\infty\cap
  \W21)(\reali; \reali)$ be a minimizer of $J$ as defined
  in~(\ref{eq:Jped}), with $S$ being the semigroup generated
  by~(\ref{eq:SCL})--(\ref{eq:Ped}).  Then, for all $r_{o} \in
  (\W1\infty\cap \W21)(\reali;\reali)$, $\rho_o$
  satisfies~(\ref{eq:last}).
\end{theorem}

\noindent The proof is deferred to Paragraph~\ref{subsec:ApplProofs}.

\section{Detailed Proofs}
\label{sec:Tech}

Below, we denote by $W_N$ the Wallis integral
\begin{equation}
  \label{eq:WN}
  W_N = \int_{0}^{\pi/2} (\cos \alpha)^N \d{\alpha} \,.
\end{equation}

\subsection{A Lemma on the Transport Equation}

The next lemma is similar to other results in recent literature, see
for instance~\cite{Ambrosio2008}.

\begin{lemma}
  \label{lem:Both}
  Let $T > 0$, so that $I=\left[0, T\right[$, and $w$ be such that
  \begin{equation}
    \label{eq:HypW}
    \begin{array}{rcl@{\qquad}rcl}
      w
      & \in &
      \C0 ( I \times \reali^N;\reali^N)
      &
      w(t)& \in & \C2(\reali^N;\reali^N) \quad \forall t\in I
      \\
      w
      & \in &
      \L\infty( I \times \reali^N;\reali^N)
      &\nabla_x w
      & \in &
      \L\infty( I \times \reali^N;\reali^{N\times N}) \,.
    \end{array}
  \end{equation}
  Assume that $R \in \L\infty \left( I; \L1(\reali^N;\reali)\right)
  \cap \L\infty \left(I \times \reali^N;\reali \right)$. Then, for any
  $r_o \in (\L1 \cap \L\infty) (\reali^N;\reali)$, the Cauchy problem
  \begin{equation}
    \label{eq:simple}
    \left\{
      \begin{array}{l}
        \partial_t r +\div \left( r \,w (t,x) \right) = R(t,x)
        \\
        r(0,x) = r_o(x)
      \end{array}
    \right.
  \end{equation}
  admits a unique Kru\v zkov solution $r \in \L\infty \bigl( I;
  \L1(\reali^N;\reali) \bigr)$, continuous from the right in time,
  given by
  \begin{equation}
    \label{eq:r}
    \!\!\!\!\!\!
    \begin{array}{rcl}
      r(t,x)
      & = &
      \displaystyle
      r_o \left(X(0; t,x) \right) \,
      \exp
      \left(
        -\int_0^t \div w \left(\tau,X(\tau; t,x) \right) \d{\tau}
      \right)
      \\
      & &
      \displaystyle
      +
      \int_0^t R \left(\tau, X(\tau; t,x) \right)
      \exp \left(-\int_\tau^t \div w\left(u, X(u; t,x) \right) \d{u} \right)
      \d{\tau} \,,
    \end{array}
    \!\!\!\!\!\!\!\!\!\!\!\!\!\!\!\!\!\!\!\!\!
  \end{equation}
  where $t \mapsto X(t;t_o,x_o)$ is the solution to the Cauchy problem
  \begin{equation}
    \label{eq:EDO}
    \left\{\begin{array}{l}
        \displaystyle
        \frac{d\chi}{dt} = w (t,\chi)
        \\[5pt]
        \chi(t_o) = x_o \,.
      \end{array}
    \right.
  \end{equation}
\end{lemma}

Note that the expression~(\ref{eq:r}) is formally justified
integrating~(\ref{eq:simple}) along the characteristics~(\ref{eq:EDO})
and obtaining
\begin{displaymath}
  \frac{d}{dt} \left(r \left(t,\chi(t) \right) \right)
  +
  r \left(t,\chi(t) \right) \div w\left(t,\chi(t) \right)
  =
  R\left(t,\chi(t) \right) .
\end{displaymath}
Recall for later use that the flow $X = X (t; t_o, x_o)$ generated
by~(\ref{eq:EDO}) can be used to introduce the change of variable $y =
X (0; t,x)$, so that $x = X (t; 0,y)$, due to standard properties of
the Cauchy problem~(\ref{eq:EDO}). Denote by $J(t,y) = \det \left(
  \nabla_y X (t; 0,y) \right)$ the Jacobian of this change of
variables. Then, $J$ satisfies the equation
\begin{equation}
  \label{eq:dJ}
  \frac{d J(t,y)}{dt} = \div w \left(t, X (t; 0,y) \right) \,J(t,y)
\end{equation}
with initial condition $J(0,y) = 1$. Hence $J(t,y) = \exp
\left(\int_0^t \div w \left(\tau, X (\tau; 0,y) \right)
  \d{\tau}\right)$ which, in particular, implies $J (t,y) > 0$ for all
$t \in I$, $y \in \reali^N$.

The natural modification to the present case
of~\cite[Definition~1]{Kruzkov} is the following.

\begin{definition}
  \label{def:SolK}
  Let $T>0$, so that $I=\left[0, T\right[$, and fix the maps $w \in
  \C0 ( I \times \reali^N; \reali)$ as in~(\ref{eq:HypW}) and $R \in
  \L\infty \bigl( I; \L1( \reali^N; \reali) \bigr) \cap \L\infty ( I
  \times \reali^N; \reali)$. Choose an initial datum $r_o \in \L\infty
  (\reali^N; \reali)$. A bounded mesurable map $r \in \L\infty \bigl(
  I ; \Lloc1(\reali^N; \reali) \bigr)$, continuous from the right in
  time, is a \emph{Kru\v zkov solution} to~(\ref{eq:simple}) if for
  all $k \in \reali$, for all test function $\phi \in \Cc\infty (]0,T[
  \times \reali^N;\rpic)$
  \begin{equation}
    \label{eq:defKruzkov2}
    \int_0^{+\infty} \!\!\! \int_{\reali^N} \!\!
    \left[
      (r - k)
      ( \pt_t \phi + w \cdot \nabla_x \phi )
      +
      ( R - k \div w ) \phi
    \right]
    \sgn (r-k) \dx \, \dt
    \geq 0
  \end{equation}
  and there exists a set $\mathcal{E}$ of zero measure in $\rpic$ such
  that for $t\in \rpic\backslash \mathcal{E}$ the function $r$ is
  defined almost everywhere in $\reali^N$ and for any $\delta>0$
  \begin{equation}
    \label{eq:Initial}
    \lim_{t\to 0, t\in ]0,T[\backslash \mathcal{E}}
    \int_{B(0,\delta)}\modulo{r(t,x)-r_o(x)}\mathrm{d}x =0\,.
  \end{equation}
\end{definition}

\begin{proofof}{Lemma~\ref{lem:Both}}
  The proof consists of several steps.

  \smallskip
  \noindent\textbf{1.}~(\ref{eq:defKruzkov2}) holds.\\
  Let $k \in \reali$ and $\phi \in \Cc\infty (\pint{I} \times
  \reali^N; \rpic)$. Then, according to Definition~\ref{def:SolK}, we
  prove~(\ref{eq:defKruzkov2}) for $r$ given as in~(\ref{eq:r}).
  By~(\ref{eq:dJ}), the semigroup property of $X$ and denoting
  $\mathcal{R}(t,y) = \int_0^t R \left( \tau, X (\tau; 0,y) \right)
  J(\tau,y) \d{\tau}$, we get
  \begin{eqnarray*}
    & &
    \int_0^{+\infty} \!\!\! \int_{\reali^N} \!\!
    \left[
      (r - k)
      ( \pt_t \phi + w \cdot \nabla_x \phi )
      +
      ( R - k \div w ) \phi
    \right]
    \sgn (r-k) \dx \, \dt
    \\[5pt]
    & = &
    \int_0^{+\infty} \!\!\! \int_{\reali^N} \!\!
    \left[
      \left(
        \frac{r_o(y)}{J(t,y)} 
        + 
        \frac{\mathcal{R}(t,y)}{J(t,y)} - k
      \right)   
    \right.
    \\
    & &
    \times
    \left( \pt_t \phi \left(t,X (t;0,y) \right)
      +
      w\left( t, X (t; 0,y) \right) \cdot
      \nabla_x \phi \left(t,X (t; 0,y) \right)
    \right)
    \\
    & &
    \left.
      +
      \left(
        R\left(t,X (t;0,y)\right)
        -
        k \div \left(  w\left(t,X (t;0, y)\right)  \right)
      \right)
      \phi \left(t,X (t;0,y)\right)
    \right]
    \\
    & &
    \qquad
    \times \sgn
    \left(
      \frac{r_o(y)}{J(t,y)} + \frac{\mathcal{R}(t,y)}{J(t,y)} - k
    \right)
    J(t,y)
    \d{y}\, \dt
    \\[5pt]
    & = &
    \int_0^{+\infty}\!\!\! \int_{\reali^N} \!\!
    \left[
      r_o(y) \frac{d}{dt} \phi \left(t,X (t;0,y) \right)
      -
      k \, J(t,y) \frac{d}{dt} \phi \left(t,X (t;0,y) \right) 
    \right.
    \\
    & &
    \left.
      -
      k \, \phi\left(t,X (t;0,y)\right) \frac{d}{dt} J(t,y)
      +
      \frac{d}{dt}\left( \mathcal{R}(t,y) \phi\left(t, X (t;0,y) \right)\right)
    \right]
    \\
    & &
    \qquad
    \times \sgn
    \left( r_o(y) + \mathcal{R} (t,y) - k\, J(t,y) \right)
    \d{y} \, \dt
    \\[5pt]
    & = &
    \int_0^{+\infty}\!\!\! \int_{\reali^N}
    \frac{d}{dt}
    \left(
      \left( r_o(y) + \mathcal{R}(t,y) - k\, J(t,y) \right) \,
      \phi \left(t,X (t;0,y) \right)
    \right)
    \\
    & &
    \times
    \sgn \left( r_o(y) + \mathcal{R}(t,y) - k\, J(t,y) \right)
    \d{y} \, \dt
    \\[5pt]
    & = &
    \int_0^{+\infty}\!\!\! \int_{\reali^N}
    \frac{d}{dt} 
    \left( 
      \modulo{r_o(y) + \mathcal{R}(t,y) - k\, J(t,y)} \, 
      \phi \left(t,X (t;0,y) \right) 
    \right)
    \d{y} \, \dt
    \\[5pt]
    & \geq &
    0 \,.
  \end{eqnarray*}

  \smallskip
  \noindent\textbf{2.}  $r \in \L\infty (I\times
  \reali^N;\reali)$.\\
  Indeed, by~(\ref{eq:r}) we have
  \begin{equation}
    \label{eq:rinfty}
    \norma{r}_{\L\infty(I\times\reali^N; \reali)}
    \leq 
    \left(
      \norma{r_o}_{\L\infty(\reali^N;\reali)}
      +
      T \norma{R}_{\L\infty(I\times\reali^N;\reali)}
    \right)
    e^{T \norma{\div w}_{\L\infty(I\times\reali^N;\reali)}} \,.
  \end{equation}

  \smallskip
  \noindent\textbf{3.} $r$ is right continuous.\\
  Consider first the case
  \begin{equation}
    \label{eq:r1}
    \left\{
      \begin{array}{l}
        \pt_t r +\div \left(r w(t,x)\right) 
        =
        0\,,
        \\
        r(0,x) = r_o(x)\,;
      \end{array}
    \right.
  \end{equation}
  where we can apply Kru\v zkov Uniqueness
  Theorem~\cite[Theorem~2]{Kruzkov}. Therefore, it is sufficient to
  show that~(\ref{eq:r}) does indeed give a Kru\v zkov solution. To
  this aim, it is now sufficient to check the continuity from the
  right at $t = 0$. Since $r_o\in (\L1 \cap \L\infty)
  (\reali^N;\reali)$, there exists a sequence $(r_{o,n})$ in
  $\Cc1(\reali^N;\reali)$ converging to $r_o$ in $\L1$.  Then, the
  corresponding sequence of solutions $(r_n)$ converges uniformly in
  time to $r$ as given by~(\ref{eq:r}). Indeed, by the same change of
  variable used above, we get
  \begin{eqnarray*}
    \int_{\reali^N} \modulo{r_n(t,x) - r(t,x)} \dx 
    & = &
    \int_{\reali^N} 
    \modulo{\frac{r_{o,n}(y)}{J(t,y)} - \frac{r_o(y)}{J(t,y)}} J(t,y) \d{y}
    \\
    & = &
    \norma{r_{o,n} - r_o}_{\L1}\,.
  \end{eqnarray*}
  Furthermore, by~(\ref{eq:r}), $r_n$ is continuous in time, in
  particular at time $t=0$. Finally, for any $\delta > 0$,
  \begin{eqnarray*}
    \int_{B(0,\delta)} \modulo{r(t,x)-r_o(x)} \dx 
    & \leq & 
    \int_{B(0,\delta)} \modulo{r(t,x)-r_n(t,x)} \dx 
    +
    \int_{B(0,\delta)} \modulo{r_n(t,x)-r_{o,n}(x)} \dx 
    \\
    & & 
    +
    \int_{B(0,\delta)} \modulo{r_{o,n}(x)-r_o(x)} \dx 
    \\
    & \leq & 
    \frac{\epsilon}{2} 
    +
    \int_{B(0,\delta)} \modulo{r_n(t,x)-r_{o,n}(x)} \dx \,, 
    \mbox{ for n large enough}
    \\
    & \leq & 
    \epsilon \,,\textrm{ for t small enough.}
  \end{eqnarray*}
  Next, we consider the system
  \begin{equation}
    \label{eq:r2}
    \left\{
      \begin{array}{l}
        \pt_t r +\div \left(r w(t,x) \right) 
        = R \,,
        \\
        r(0,x) = 0 \,,
      \end{array}
    \right.
  \end{equation}
  and introduce the map
  \begin{displaymath}
    \mathcal{F}(t,h,\tau,y)
    =
    \exp \left[
      \int_0^t \div w \left(u, X(u;0,y) \right) \d{u}
      -
      \int_\tau^{t+h} 
      \!\!\!
      \div w \left(u, X\left(u;t+h,X(t;0,y) \right)\right)
      \d{u}
    \right]
  \end{displaymath}
  so that the solution to~(\ref{eq:r2}) satisfies by~(\ref{eq:r}), for
  all $h>0$,
  \begin{eqnarray*}
    & &
    \norma{r(t+h)-r(t)}_{\L1}
    \\
    & = &
    \int_{\reali^N} 
    \left|
      \int_0^{t+h} R \left(\tau, X \left(\tau;t+h,X(t;0,y) \right) \right)
      \mathcal{F}(t,h,\tau,y)
      \d{\tau} \right.
    \\
    &&\left.\qquad-
      \int_0^t R\left(\tau,X(\tau;0,y) \right) J(\tau,y) \d{\tau}
    \right| \d{y}
    \\
    & \leq &
    \int_{\reali^N} \int_0^t 
    \modulo{R \left( \tau, X \left( \tau;t+h,X(t;0,y)\right)\right)} 
    \modulo{\mathcal{F}(t,h,\tau,y) - J(\tau,y)}
    \d\tau \d{y}
    \\
    & & 
    +
    \int_{\reali^N} \int_0^t 
    J(\tau,y) 
    \modulo{
      R\left(\tau, X \left( \tau;t+h,X(t;0,y)\right)\right)
      -
      R \left( \tau,X(\tau;0,y)\right)
    } 
    \d\tau \d{y}
    \\
    & &
    +
    \int_t^{t+h} \int_{\reali^N} 
    \mathcal{F}(t,h,\tau,y) 
    \modulo{R \left(\tau,X \left( \tau;t+h,X(t;0,y)\right)\right)}
    \d{y} \d\tau  \,.
  \end{eqnarray*}
  The former summand above vanishes as $h \to 0$ because the integrand
  is uniformly bounded in $\L1$ and converges pointwise to $0$, since
  $X\left(u;t+h,X(t;0,y) \right) \stackrel{h\to 0}{\longrightarrow}
  X(u; 0,y)$ and also $\mathcal{F}(t,h,\tau,y) \stackrel{h\to
    0}{\longrightarrow} J(\tau,y)$.  The second one, in the same
  limit, vanishes by the Dominated Convergence Theorem, $R$ being in
  $\L1$ and by the boundedness of $J$. Indeed, if $(R_n)$ is a
  sequence of functions in $\Cc1(\reali^N;\reali)$ that converges to
  $R$ in $\L1$ we have
  \begin{eqnarray*}
    & &\int_{\reali^N} \int_0^t 
    J(\tau,y) 
    \modulo{
      R\left(\tau, X \left( \tau;t+h,X(t;0,y)\right)\right)
      -
      R \left( \tau,X(\tau;0,y)\right)
    } 
    \d\tau \d{y} \\
    & \leq &   \int_{\reali^N} \int_0^t 
    J(\tau,y) 
    \modulo{
      R_n\left(\tau, X \left( \tau;t+h,X(t;0,y)\right)\right)
      -
      R_n \left( \tau,X(\tau;0,y)\right)
    } 
    \d\tau \d{y}\\
    & &+  \int_{\reali^N} \int_0^t 
    J(\tau,y) 
    \modulo{
      R_n\left(\tau, X \left( \tau;t+h,X(t;0,y)\right)\right)
      -
      R \left( \tau,X(\tau;t+h, X(t;0,y)\right)
    } 
    \d\tau \d{y}\\
    & & +  \int_{\reali^N} \int_0^t 
    J(\tau,y) 
    \modulo{
      R\left(\tau, X \left( \tau;0,y)\right)\right)
      -
      R_n \left( \tau,X(\tau;0,y)\right)
    } 
    \d\tau \d{y}\,.
  \end{eqnarray*}
  $J(\tau,y)$ is uniformly bounded on $[0,t]\times \reali^N$. We can
  first fix $n$ large enough so that the second and third terms will
  be small, independently of $h$. Then, taking $h$ small enough, we
  know from Dominated Convergence Theorem that the first term will
  shrink to 0.  The integrand in the latter summand is in $\L\infty$
  since $R$ is in $\L1$.

  In general, right continuity follows by linearity adding the
  solutions to~(\ref{eq:r1}) and~(\ref{eq:r2}).

  \smallskip
  \noindent\textbf{4.}  $r \in \L\infty
  \bigl(I; \L1(\reali^N; \reali) \bigr)$.\\
  Indeed, for all $t \in I$ we have
  \begin{displaymath}
    \norma{r(t)}_{\L1}
    \leq 
    \left(
      \norma{r_o}_{\L1} 
      +
      t \norma{R}_{\L\infty(I;\L1)}
    \right)
    \exp \left( t \norma{\div w}_{\L\infty} \right) \,.
  \end{displaymath}

  \smallskip
  \noindent\textbf{5.} The solution to~(\ref{eq:simple}) is unique.

  First, assume that $w \in \C2 (I\times \reali^N; \reali^N)$ and $R
  \in \C2 (I \times \reali^N; \reali)$. Then, Kru\v zkov Uniqueness
  Theorem~\cite[Theorem~2]{Kruzkov} applies.

  Second, assume that $w \in \C2 ( I \times \reali^N; \reali^N)$ and
  $R$ satisfies the present assumptions. Then, we use the same
  procedure as in the proof
  of~\cite[Theorem~2.6]{ColomboMercierRosini}. There, the general
  scalar balance law $\partial_t u + \div f(t,x,u) = F(t,x,u)$ is
  considered, under assumptions that allow first to apply Kru\v zkov
  general result and, secondly, to prove stability estimates on the
  solutions. Remark that these latter estimates are proved therein
  under the only requirement that solutions are Kru\v zkov solutions,
  according to~\cite[Definition~1]{Kruzkov} or, equivalently,
  Definition~\ref{def:SolK}. Here, the existence part has been proved
  independently from Kru\v zkov result and under weaker assumptions.

  Let $(R_n)$ be a sequence in $\Cc2$ that converges in $\L1$ to $R
  \in \C0 \bigl( I; \L1(\reali^N,\reali) \bigr)$. Also with reference
  to the notation of~\cite[Theorem~2.6]{ColomboMercierRosini},
  consider the equations
  \begin{eqnarray}
    \label{eq:fF}
    \partial_t r_n + \div \left( r_n \, w(t,x) \right) = R_n(t,x)
    & \mbox{and let}&
    \left\{
      \begin{array}{rcl}
        f(t,x,r) & = & r \, w(t,x)
        \\
        F(t,x,r) & = & R_n(t,x)
      \end{array}
    \right.
    \\
    \label{eq:gG}
    \partial_t r + \div \left( r \, w(t,x) \right) = R(t,x)
    & \mbox{and let}&
    \left\{
      \begin{array}{rcl}
        g(t,x,r) & = &  r \, w(t,x)
        \\
        G(t,x,r) & = & R(t,x)
      \end{array}
    \right.
  \end{eqnarray}
  with the same initial datum $r_o \in (\L1 \cap \L\infty)
  (\reali^N;\reali)$.

  Note that here the sources $F$ and $G$ do not depend on $r$, hence
  the proof of~\cite[Theorem~2.6]{ColomboMercierRosini} can be
  repeated with $G \in \C0 \bigl( I; \L1(\reali^N,\reali) \bigr)$
  instead of $\C0 \bigl( I \times \reali^N ; \reali \bigr)$. Indeed,
  in the proof of~\cite[Theorem~2.6]{ColomboMercierRosini}, it is
  sufficient to have the $(t,x)$--regularity in the source term $F$ of
  the first equation and existence and continuity of the derivative
  $\pt_r (F - G)$, which here vanishes. Besides, here the two flows
  $f$ and $g$ are identical, hence we do not need the $\BV$ estimate
  provided by~\cite[Theorem~2.5]{ColomboMercierRosini}.

  Thus, to apply the stability estimate
  in~\cite[Theorem~2.6]{ColomboMercierRosini}, we are left to check
  the following points:
  \begin{itemize}
  \item the derivatives $\pt_r f = w$, $\pt_r \nabla_x f = \nabla_x
    w$, $\nabla^2_x f = r \nabla^2_x w$, $\pt_r F$ and $\nabla_x F$
    exist and are continuous;
  \item $\pt_r f = w$ and $F - \div f = F - r \div w$ are bounded in
    $I \times \reali^N \times[-A,A]$ for all $A \geq 0$;
  \item $\pt_r (F - \div f)(t,x,r) = -\div w \in \L\infty (I \times
    \reali^N\times \reali ; \reali)$,
  \item $\nabla_x\pt_r f(t,x,r) = \nabla_x w \in \L\infty (I \times
    \reali^N\times\reali; \reali^{N\times N})$.
  \end{itemize}
  \noindent Hence, if $r$ is any solution to~(\ref{eq:gG}) and $r_n$
  is the solution to~(\ref{eq:fF}) in the sense of
  Definition~\ref{def:SolK}, then for $t \in I$, $x_o \in \reali^N$,
  $\delta \geq 0$, $M = \norma{w}_{\L\infty (\rpic \times \reali^N;
    \reali^N)}$:
  \begin{displaymath}
    \norma{(r_n - r)(t)}_{\L1(B(x_o,\delta); \reali)}
    \leq 
    \int_0^t e^{\kappa t} 
    \norma{(R_n - R)(s)}_{\L1(B(x_o,\delta+M(t-s));\reali)} \d{s}\, ,
  \end{displaymath}
  where $\kappa = 2 N \norma{\nabla_x w}_{\L\infty(I\times \reali^N;
    \reali^{N\times N})}$. Therefore, if $(R_n)$ converges in $\L1$ to
  $R$, then $(r_n)$ is a Cauchy sequence in $\L\infty \bigl( I;
  \L1(\reali^N; \reali) \bigr)$ and $r \in \L \infty \bigl( I;
  \L1(\reali^N;\reali) \bigr)$ is uniquely characterized as its limit.

  Third, we consider the general case. Again, we rely on the proof
  of~\cite[Theorem~2.6]{ColomboMercierRosini} extending it to the case
  of $w \in \C0 (I \times \reali^N;\reali^N)$. Indeed, therein the
  higher regularity in time of the flow is used to apply Kru\v zkov
  Existence Theorem~\cite[Theorem~5]{Kruzkov}, to prove the $\BV$
  estimates in~\cite[Theorem~2.4]{ColomboMercierRosini} and to obtain
  the limit~\cite[(5.11)]{ColomboMercierRosini}. In the former case,
  our existence proof in the previous steps replaces the use of Kru\v
  zkov result. $\BV$ estimates are here not necessary, for we keep
  here the flow fixed.  In the latter case, a simple argument based on
  the Dominated Convergence Theorem allows to get the same limit.
\end{proofof}

Remark that as an immediate corollary of Lemma~\ref{lem:Both} we obain
that any solution to~(\ref{eq:simple}) in the sense of
Definition~\ref{def:SolK} is represented by~(\ref{eq:r}).

\subsection{Proof of Theorem~\ref{thm:main}}

\begin{lemma}
  \label{lem:Infty}
  Let $T > 0$, so that $I=\left[0, T\right[$, and $w$ be as
  in~(\ref{eq:HypW}) such that
  \begin{eqnarray}
    \label{eq:HypW2}
    \div w
    & \in &
    \L\infty \left(I; \L1(\reali^N;\reali) \right)
    \\
    \label{eq:HypW3}
    \nabla_x \div w
    & \in &
    \L\infty \left( I; \L1(\reali^N;\reali^N) \right) \,.
  \end{eqnarray}
  Then, for any $\rho_o \in (\L1 \cap \L\infty) (\reali^N; \reali)$,
  the Cauchy problem~(\ref{eq:simple}) with $R=0$ admits a unique
  solution $\rho \in \L\infty \left( I; \L1 (\reali^N; \reali)
  \right)$, right continuous in time and satisfying
  \begin{equation}
    \label{eq:BoundInfty}
    \norma{\rho(t)}_{\L\infty(\reali^N;\reali)} 
    \leq
    \norma{\rho_o}_{\L\infty(\reali^N;\reali)} \,
    \exp \left(
      t \, \norma{\div w}_{\L\infty([0,t]\times \reali^N;\reali)}
    \right)
  \end{equation}
  for all $t \in I$.  Moreover, this solution has the following
  properties:
  \begin{enumerate}
  \item
    \label{it:Infty:1}
    \begin{tabular}[t]{lcp{0.7\linewidth}}
      $\rho_o \geq 0$ a.e.
      & $\Rightarrow$
      & $\rho(t) \geq 0$ a.e., for all $t \in I$
      \\[6pt]
      $\div w\geq 0$ a.e.
      & $\Rightarrow$
      & $\displaystyle \norma{\rho(t)}_{\L\infty}\leq
      \norma{\rho_o}_{\L\infty}$, for all $t \in I$.
    \end{tabular}
  \item \label{it:Infty:2} If $\rho_o \in \mathcal{X}$ then, for all
    $t \in I$, we have $\rho(t) \in \mathcal{X}$ and setting $\kappa_o
    = N W_N (2N+1) \norma{\nabla_x
      w}_{\L\infty(I\times\reali^N;\reali^{N\times N})}$, we also get
    \begin{eqnarray*}
      \tv\left( \rho(t) \right) 
      & \leq &
      \tv (\rho_o) e^{\kappa_o t} 
      \\
      & &
      +
      N W_N \int_0^t e^{\kappa_o(t-s)} 
      \int_{\reali^N} e^{s \norma{\div w}_{\L\infty}}
      \norma{\nabla_x \div w(s,x)} \, \dx \, \d{s} \,
      \norma{\rho_o}_{\L\infty} \,.
    \end{eqnarray*}
    Furthermore, $\rho \in \C0 \bigl( I; \L1 (\reali^N; \reali)
    \bigr)$.
  \item \label{it:Infty:3} If $\rho_1$, $\rho_2$ are the solutions
    of~(\ref{eq:simple}) associated to $w_1$, $w_2$ with $R_1 = R_2 =
    0$ and with initial conditions $\rho_{1,o}$, $\rho_{2,o}$ in
    $\mathcal{X}$, then for all $t \in I$
    \begin{eqnarray*}
      & &
      \norma{(\rho_1-\rho_2)(t)}_{\L1}
      \\
      & \leq &
      e^{\kappa t}\norma{\rho_{1,o}-\rho_{2,o}}_{\L1}
      +
      \frac{e^{\kappa_o t}-e^{\kappa t}}{\kappa_o-\kappa}
      \tv(\rho_{1,o}) 
      \norma{w_1-w_2}_{\L\infty}
      \\
      & &
      + 
      N W_N \int_0^t 
      \frac{e^{\kappa_o (t-s)} - e^{\kappa(t-s)}}{\kappa_o-\kappa}
      \int_{\reali^N} e^{s \norma{\div w}_{\L\infty}}
      \norma{\nabla_x\div w_1 (s,x)}\, \dx \, \d{s} 
      \\
      & &
      \qquad\quad 
      \times 
      \norma{\rho_o}_{\L\infty} \,
      \norma{w_1-w_2}_{\L\infty} 
      \\
      & &
      + 
      \int_0^t e^{\kappa(t-s)} e^{s \norma{\div w}_{\L\infty}}
      \int_{\reali^N} \modulo{\div(w_1-w_2)(s,x)} \, \dx \, \d{s}
      \norma{\rho_o}_{\L\infty} \,,
    \end{eqnarray*}
    where $\kappa = 2 N \norma{\nabla_x w_1}_{\L\infty(I \times
      \reali^N; \reali^{N\times N})}$ and $\kappa_o$ as
    in~\ref{it:Infty:2} above.
  \item \label{it:Infty:4} If there exists $C \geq 0$ such that
    \begin{equation}
      \label{eq:CondC}
      \norma{\nabla_x^2 w}_{\L\infty(I\times\reali^N; \reali^{N\times N\times N})}
      \leq
      C
      \quad \mbox{ and } \quad
      \norma{\nabla_x w}_{\L\infty(I\times\reali^N; \reali^{N\times N})} 
      \leq
      C
    \end{equation}
    then
    \begin{displaymath}
      \begin{array}{lcl}
        \rho_o \in \W11(\reali^N; \reali)
        & \Rightarrow &
        \left\{
          \begin{array}{l}
            \rho(t) \in \W11(\reali^N; \reali)
            \quad \mbox{ for all } t \in I
            \\[4pt]
            \norma{\rho(t)}_{\W11}
            \leq 
            e^{2Ct}
            \norma{\rho_o}_{\W11} \,,
          \end{array}
        \right.
        \\[16pt]
        \rho_o
        \in \W1\infty(\reali^N; \reali)
        & \Rightarrow &
        \left\{
          \begin{array}{l}
            \rho(t)\in\W1\infty(\reali^N; \reali)
            \quad \mbox{ for all } t \in I
            \\[4pt]
            \norma{\rho(t)}_{\W1\infty}
            \leq 
            e^{2Ct} \, \norma{\rho_o}_{\W1\infty} \, .
          \end{array}
        \right.
      \end{array}
    \end{displaymath}
  \item \label{it:Infty:5} If there exists $C\geq 0$ such
    that~(\ref{eq:CondC}) holds together with $\norma{\nabla_x^3
      w}_{\L\infty(I\times\reali^N; \reali^{N\times N\times N\times
        N})} \leq C$, then $\rho_o\in \W21(\reali^N; \reali)$ implies
    \begin{displaymath}
      \rho(t) \in \W21(\reali^N; \reali) 
      \quad \mbox{ for all } t \in I
      \quad \mbox{ and } \quad
      \norma{\rho(t)}_{\W21}
      \leq
      (1+Ct)^2 e^{3Ct} \,    
      \norma{\rho_o}_{\W21} \,.
    \end{displaymath}
  \end{enumerate}
\end{lemma}

\begin{proof}
  The existence of a Kru\v zkov solution follows from
  Lemma~\ref{lem:Both}.  But we can also refer to~\cite[theorems~2
  and~5]{Kruzkov}, the assumptions in~\cite[\S~5]{Kruzkov} being
  satisfied thanks to~(\ref{eq:HypW}).  The $\L\infty$ bound directly
  follows from~(\ref{eq:r}), which now reads
  \begin{equation}\label{eq:rho}
    \rho(t,x) 
    = 
    \rho_o \left(X(0;t,x) \right)  \, 
    \exp \left(
      -\int_0^t \div w \left(\tau,X(\tau; t,x) \right) \d{\tau} 
    \right)\,.
  \end{equation}

  The representation formula~(\ref{eq:rho}) also implies the bounds
  at~\ref{it:Infty:1}.

  The bound on the total variation at~\ref{it:Infty:2} follows
  from~\cite[Theorem~2.5]{ColomboMercierRosini}, the hypotheses on $w$
  being satisfied thanks to~(\ref{eq:HypW}) and~(\ref{eq:HypW3}). More
  precisely, we do not have here the $\C2$ regularity in time as
  required in~\cite[Theorem 2.5]{ColomboMercierRosini}, but going
  through the proof of this result, we can see that only the
  continuity in time of the flow function $f(t,x,r)= r w(t,x)$ is
  necessary. Indeed, time derivatives of $f$ appear in the proof
  of~\cite[Theorem 2.5]{ColomboMercierRosini} when we bound the terms
  $J_t$ and $L_t$, see~\cite[between~(4.18)
  and~(4.19)]{ColomboMercierRosini}. However, the use of the Dominated
  Convergence Theorem allows to prove that $J_t$ and $L_t$ converge to
  zero when $\eta$ goes to $0$ without any use of time derivatives.
  The continuity in times follows
  from~\cite[Remark~2.4]{ColomboMercierRosini}, thanks
  to~(\ref{eq:HypW2}) of $w$ and the bound on the total variation.

  Similarly, the stability estimate at~\ref{it:Infty:3} is based
  on~\cite[Theorem~2.6]{ColomboMercierRosini}. Indeed, we use once
  again a flow that is only $\C0$ instead of $\C2$ in time. Besides,
  in the proof of~\cite[Theorem~2.6]{ColomboMercierRosini}, the
  $\L\infty$ bound into the integral term
  in~\cite[Theorem~2.6]{ColomboMercierRosini} can be taken only in
  space, keeping time fixed. With this provision, the proof
  of~\ref{it:Infty:3} is exactly the same as that
  in~\cite{ColomboMercierRosini}, so we do not reproduce it here. The
  same estimate is thus obtained, except that the $\L\infty$ bound of
  the integral term is taken only in space.

  The proofs of the $\W11$ and $\W1\infty$ bounds
  at~{\ref{it:Infty:4}} are similar. They follow from the
  representation~(\ref{eq:rho}), noting that $\norma{\nabla_x
    X}_{\L\infty} \leq e^{Ct}$. Indeed,
  \begin{eqnarray*}
    \nabla_x
    X(t;0,x) 
    & = &
    \mathbf{Id} 
    + 
    \int_0^t \nabla_x
    w(\tau;X(\tau;0,x)) \nabla_xX(\tau;0,x) \d{\tau}
    \,,\mbox{ hence}
    \\
    \norma{ \nabla_x X(t;0,x)}
    & \leq  &   %% changed: =
    1 
    +
    \int_0^t \norma{\nabla_x w(\tau;X(\tau;0,x))} 
    \norma{\nabla_x X(\tau;0,x)} \d{\tau}
    \\
    & \leq &
    1
    +
    \int_0^tC\norma{\nabla_x X(\tau;0,x)} \d{\tau}
  \end{eqnarray*}
  and a direct application of Gronwall Lemma gives the desired
  bound. Hence, we obtain
  \begin{displaymath}
    \norma{\nabla\rho(t)}_{\L\infty} 
    \leq 
    (e^{2Ct}-e^{Ct}) \, \norma{\rho_o}_{\L\infty}
    +
    e^{2Ct} \, \norma{\nabla\rho_o}_{\L\infty}
  \end{displaymath}
  and consequently
  \begin{displaymath}
    \norma{\rho(t)}_{\W1\infty}
    \leq
    e^{2Ct} \, \norma{\rho_o}_{\W1\infty}\,.
  \end{displaymath}
  The $\L1$ estimate is entirely analogous.

  The $\W21$ bound at~\ref{it:Infty:5}.~also comes from
  the~(\ref{eq:rho}). Indeed, again thanks to Gronwall Lemma, we get
  $\norma{\nabla^2_x X}_{\L\infty} \leq e^{2Ct}-e^{Ct}$. Using the
  estimates above, together with
  \begin{displaymath}
    \norma{\nabla^2 \rho(t)}_{\L1}
    \leq
    (2e^{2Ct} -3e^{Ct}+1) e^{Ct} \, \norma{\rho_o}_{\L1}
    +
    3(e^{Ct}-1) e^{2Ct} \, \norma{\nabla\rho_o}_{\L1}
    + 
    e^{3Ct} \, \norma{\nabla^2 \rho_o}_{\L1}\,,
  \end{displaymath}
  we obtain
  \begin{displaymath}
    \norma{\rho(t)}_{\W21}
    \leq
    (2e^{Ct}-1)^2 e^{Ct} \,    
    \norma{\rho_o}_{\W21} 
  \end{displaymath}  
  concluding the proof.
\end{proof}

We use now these tools in order to obtain the existence of a solution
for~(\ref{eq:SCL}).

\begin{proofof}{Theorem~\ref{thm:main}}
  Fix $\alpha, \beta>0$ with $\beta>\alpha$. Let $T_* =
  \left(\ln(\beta/\alpha) \right)/C(\beta)$, with $C$ as
  in~\textbf{(V1)}.  Define the map
  \begin{displaymath}
    \mathcal{Q}\;\;
    \colon \;\;
    \begin{array}{ccc}
      \C0 \left( I_*; \mathcal{X}_\beta \right)
      & \to &
      \C0 \left( I_*; \mathcal{X}_\beta \right)
      \\
      \sigma
      & \mapsto &
      \rho
    \end{array}
  \end{displaymath}
  where $I_* = \left[0, T_*\right[$ and $\rho$ is the Kru\v zkov
  solution to
  \begin{equation}
    \label{eq:w}
    \left\{
      \begin{array}{l}
        \partial_t \rho + \div ( \rho \, w ) = 0
        \\
        \rho(0,x) = \rho_o(x)\,
      \end{array}
    \right.
    \qquad \mbox{ with } \quad
    \begin{array}{rcl}
      w & = & V(\sigma) 
      \\
      \rho_o & \in & \mathcal{X}_\alpha\,.
    \end{array}
  \end{equation}
  The assumptions~\textbf{(V1)} imply the hypotheses on $w$ necessary
  in Lemma~\ref{lem:Infty}. Therefore, a solution $\rho$
  to~(\ref{eq:w}) exists and is unique. In particular, the continuity
  in time of $\rho$ follows from~\ref{it:Infty:2}~in
  Lemma~\ref{lem:Infty}, due to the boundedness of the total
  variation. Note that by~(\ref{eq:BoundInfty}), the choice of $T_*$
  and~\textbf{(V1)}, $\norma{\rho(t)}_{\L\infty} \leq \beta$ and hence
  $\mathcal{Q}$ is well defined.

  Fix $\sigma_1,\sigma_2$ in $\C0 (I_*; \mathcal{X}_\beta)$. Call $w_i
  = V(\sigma_i)$ and $\rho_i$ the corresponding solutions. With the
  same notations of~\cite[Theorem~2.6]{ColomboMercierRosini}, we let
  \begin{displaymath}
    \kappa_o
    =
    N \, W_N(2N+1) \, 
    \norma{\nabla_x w_1}_{\L\infty(I_*\times \reali^N;\reali^{N\times N})}\, ,
    \qquad
    \kappa
    =
    2 \, N \, 
    \norma{\nabla_x w_1}_{\L\infty(I_*\times \reali^N;\reali^{N\times N})}\, .
  \end{displaymath}
  Note that by~(\ref{eq:WN})
  \begin{displaymath}
    \frac{\kappa_o}{\kappa}
    % = \left( N + \frac{1}{2} \right) W_N
    \geq
    \left( N + \frac{1}{2} \right)
    \int_0^{\pi/2} \left( 1-\frac{2}{\pi}x \right)^N
    \,
    \dx
    =
    \frac{\pi}{2} \left( 1- \frac{1}{2(N+1)}\right)
    \geq
    \frac{3\pi}{8} >
    1
  \end{displaymath}
  hence $\kappa_o > \kappa$. Then, by 4 of Lemma~\ref{lem:Infty}
  and~\textbf{(V1)}, we obtain a bound on $\kappa_o$. Indeed,
  \begin{displaymath}
    \norma{\nabla_x V(\sigma_1)}_{\L\infty(I_*\times\reali^N;\reali^{N\times N})}
    \leq
    C \left(
      \norma{\sigma_1}_{\L\infty(I_*\times\reali^N;\reali)} 
    \right),
  \end{displaymath}
  and since $\sigma_1 \in \C0 ( I_*; \mathcal{X}_\beta )$, finally
  $\kappa_o \leq N W_N ( 2 N + 1) \, C(\beta)$.  Let us denote
  \begin{equation}
    \label{C'}
    C 
    =
    C (\beta)
    \quad \textrm{ and } \quad
    C' = N \, W_N \, (2N+1) \, C (\beta)
    \,.
  \end{equation}
  Again, \textbf{(V1)} implies the following uniform bounds on all
  $\sigma_1, \sigma_2 \in \C0 ( I_*; \mathcal{X}_\beta)$:
  \begin{eqnarray*}
    \norma{\nabla_x^2 V(\sigma_1)}_{\L\infty(I_*;\L1(\reali^N;\reali^{N\times N \times N}))}
    & \leq &  
    C \,,
    \\
    \norma{V(\sigma_1)-V(\sigma_2)}_{\L\infty(I_*\times\reali^N;\reali^N)}
    & \leq &
    C \, \norma{\sigma_1-\sigma_2}_{\L\infty(I_*;\L1(\reali^N;\reali))}\,,
    \\
    \norma{\div\left(V(\sigma_1) - V(\sigma_2)\right)}_{\L\infty(I_*;\L1(\reali^N;\reali))}
    & \leq &
    C \, \norma{\sigma_1-\sigma_2}_{\L\infty(I_*;\L1(\reali^N;\reali))}\,.
  \end{eqnarray*}
  Thus, we can apply~\cite[Theorem~2.6]{ColomboMercierRosini}. We get,
  for all $t \in I_*$,
  \begin{eqnarray*}
    \norma{(\rho_1-\rho_2)(t)}_{\L1}
    & \leq &
    C t e^{C't}\tv(\rho_o) \norma{\sigma_1-\sigma_2}_{\L\infty([0,t]; \L1)}
    \\
    & &
    +  
    C^2 N W_N e^{Ct} \int_0^t  (t-s)e^{C'(t-s)} \,  \d{s}
    \, \norma{\rho_{o,1}}_{\L\infty} \norma{\sigma_1-\sigma_2}_{\L\infty([0,t];\L1)} 
    \\
    & &
    +
    e^{Ct} 
    \int_0^t C \, e^{ C'(t-s)} \,
    \norma{(\sigma_1-\sigma_2)(s)}_{\L1} \,  \d{s}
    \norma{\rho_{o,1}}_{\L\infty} \, .
  \end{eqnarray*}
  Therefore, we obtain the following Lipschitz estimate:
  \begin{eqnarray*}
    & & 
    \norma{\mathcal{Q} (\sigma_1) - \mathcal{Q}(\sigma_2)}_{\L\infty(I;\L1)}
    \\
    & \leq & 
    C T e^{C'T}
    \left[ 
      \tv(\rho_o) + (N W_N C T + 1 ) e^{CT} \norma{\rho_{o,1}}_{\L\infty}
    \right] 
    \norma{\sigma_1-\sigma_2}_{\L\infty(I; \L1)}.
  \end{eqnarray*}
  Here we introduce the strictly increasing function
  \begin{displaymath}
    f(T)
    = 
    CTe^{C'T}
    \left[ 
      \tv(\rho_o)+(N W_N C T + 1 ) e^{CT} \norma{\rho_{o,1}}_{\L\infty}
    \right]
  \end{displaymath}
  and we remark that $f(T) \to 0$ when $T \to 0$. Choose now $T_1 > 0$
  so that $f(T_1) = 1/2$. Banach Contraction Principle now ensures the
  existence and uniqueness of a solution $\rho^*$ to~(\ref{eq:SCL}) on
  $[0,\bar T]$ in the sense of Definition~\ref{def:sol}, with $\bar T
  = \min\{T_*,T_1\}$. In fact, if $T_1<T_*$, we can prolongate the
  solution until time $T_*$. Indeed, if we take $\rho^*(T_1)$ as
  initial condition, we remark that $\norma{\rho^*(T_1)}_{\L\infty}
  \leq \norma{\rho_o}_{\L\infty} e^{C(\beta)T_1}$. Consequently, the
  solution of~(\ref{eq:w}) on $[T_1,T_*]$ instead of $I_*$ satisfy,
  thanks to~(\ref{eq:BoundInfty})
  \begin{displaymath}
    \norma{\rho(t)}_{\L\infty}
    \leq 
    \norma{\rho^*(T_1)}_{\L\infty} e^{C(\beta)(t-T_1)}
    \leq  
    \norma{\rho_o}_{\L\infty} e^{C(\beta)T_1}
    e^{C(\beta)(t-T_1)} 
    \leq 
    \norma{\rho_o}_{\L\infty} e^{C(\beta)T_*}\,,
  \end{displaymath}
  which is less than $\beta$ thanks to the definition of $T_{*}$ and
  since $\rho_{o} \in \mathcal{X}_{\alpha}$.

  Now, we have to show that $T_n \geq T_*$ for $n$ sufficiently
  large. To this aim, we obtain the contraction estimate
  \begin{eqnarray*}
    & &
    \norma{\mathcal{Q} (\sigma_1) - \mathcal{Q}(\sigma_2)}_{\L\infty([T_{n},T_{n+1}];\L1)}
    \\
    & \leq &
    C(T_{n+1}-T_n)e^{C'(T_{n+1}-T_n)}
    \left[ 
      \tv \left(\rho(T_n) \right) + (N W_N C ( T_{n+1} - T_n) +1)e^{CT_n}\norma{\rho_o}_{\L\infty}
    \right]
    \\
    & &
    \times 
    \norma{\sigma_1-\sigma_2}_{\L\infty([T_{n},T_{n+1}]; \L1)}
    \\
    & \leq &
    \left[ 
      \tv(\rho_o)\,e^{C'T_n}+C'T_n e^{C'T_n} +(NW_N C(T_{n+1}-T_n)
      +1)e^{CT_n}\norma{\rho_o}_{\L\infty}\right]
    \\
    & &
    \times   C'(T_{n+1}-T_n)e^{C'(T_{n+1}-T_n)}
    \norma{\sigma_1-\sigma_2}_{\L\infty([T_{n},T_{n+1}]; \L1)}
  \end{eqnarray*}
  where we used the bounds on $\tv \left( \rho(T_n) \right)$ and
  $\norma{\rho(T_n)}_{\L\infty(\reali^N;\reali)}$ provided by
  Lemma~\ref{lem:Infty} associated to the conditions~\textbf{(A)}
  and~\textbf{(V1)}. We may thus extend the solution up to time
  $T_{n+1}$, where we take $T_{n+1} > T_n$ such that
  \begin{eqnarray*}
    \left[
      \tv(\rho_o) \, e^{C'T_n} 
      +
      C T_n e^{C'T_n} 
      +
      \left( N W_N C(T_{n+1}-T_n) + 1 \right) e^{CT_n}
      \norma{\rho_o}_{\L\infty}
    \right]
    \times
    \\
    \times
    C (T_{n+1}-T_n) e^{C'(T_{n+1}-T_n)} 
    & = &
    \frac{1}{2}\,.
  \end{eqnarray*}
  If the sequence $(T_n)$ is bounded, then the left hand side above
  tends to $0$, whereas the right hand side is taken equal to $1/2 >
  0$. Hence, the sequence $(T_n)$ is unbounded. In particular, for $n$
  large enough, $T_n$ is larger than $T_*$; thus the solution
  to~(\ref{eq:SCL}) is defined on all $I_*$.

  The Lipschitz estimate follows by applying the same procedure as
  above, in the case when the initial conditions are not the same.

  The $\L\infty$ and $\tv$ bounds follow from~(\ref{eq:BoundInfty})
  and from point~\ref{it:Infty:2}~in Lemma~\ref{lem:Infty}.
\end{proofof}

The proof of Lemma~\ref{lem:Invariance} directly follows from the
second bound in~\ref{it:Infty:1}.~of Lemma~\ref{lem:Infty}.

\begin{proofof}{Theorem~\ref{thm:V0}}
  We consider the assumptions~\textbf{(A)} and~\textbf{(B)}
  separately.

  \smallskip
  \noindent\textbf{(A)}:\quad Let $T > 0$, so that $I=\left[0,
    T\right[$, and fix a positive $\alpha$. As in the proof of
  Theorem~\ref{thm:main}, we define the map
  \begin{displaymath}
    \mathcal{Q}\;\;
    \colon \;\;
    \begin{array}{ccc}

      \C0 \left( I; \mathcal{X}_\alpha \right)
      & \to &
      \C0 \left( I; \mathcal{X}_\alpha \right)
      \\
      \sigma
      & \mapsto &
      \rho
    \end{array}
  \end{displaymath}
  where $\rho$ is the Kru\v zkov solution to~(\ref{eq:w}) with
  $\rho_o\in \mathcal{X}_\alpha$. The existence of a solution
  for~(\ref{eq:w}) in $\L\infty(I,\L1(\reali^N;\reali))$ is given by
  Lemma~\ref{lem:Infty}, the set of assumptions~\textbf{(V1)} allowing
  to check the hypotheses on $w$. Note that furthermore~\textbf{(A)}
  gives an $\L\infty$ bound on $\rho$, thanks to
  Lemma~\ref{lem:Invariance}, so that for all $t\in I$, $\rho(t) \in
  [0,\alpha]$, a.e.~in $x$. Fix $\sigma_1,\sigma_2$ in $\C0 \left( I;
    \mathcal{X}_\alpha \right)$, call $w_i = V(\sigma_i)$ and let
  $\rho_1$, $\rho_2$ be the associated solutions. With the same
  notations of~\cite[Theorem~2.6]{ColomboMercierRosini}, we let as in
  the proof of Theorem~\ref{thm:main},
  \begin{displaymath}
    \kappa_o
    =
    N \, W_N(2N+1) \, 
    \norma{\nabla_x w_1}_{\L\infty(I\times \reali^N;\reali^{N\times N})}\, ,
    \qquad
    \kappa
    =
    2 \, N \, 
    \norma{\nabla_x w_1}_{\L\infty(I\times \reali^N;\reali^{N\times N})}\, .
  \end{displaymath}
  so that $\kappa_o > \kappa$. Then we use Lemma~\ref{lem:Infty} and
  assumptions~\textbf{(V1)} in order to find a bound on
  $\kappa_o$. Indeed, by~\textbf{(V1)} we have:
  \begin{displaymath}
    \norma{\nabla_x V(\sigma_1)}_{\L\infty(I\times\reali^N;\reali^{N\times N})}
    \leq
    C \left(
      \norma{\sigma_1}_{\L\infty(I\times\reali^N;\reali^{N\times N})} 
    \right),
  \end{displaymath}
  and since $\sigma_1 \in \C0 \left( I;\mathcal{X}_\alpha \right)$, we
  have $\norma{\sigma_1}_{\L\infty} \leq \alpha$ so
  that $\kappa_o \leq N W_N ( 2 N + 1) C(\alpha)$.  Denote
  \begin{equation}
    \label{C'V0}
    C' = N W_N (2N+1) C(\alpha) 
    \quad \textrm{ and } \quad
    C = C(\alpha) \,.
  \end{equation}
  The following bounds are also available uniformly for all $\sigma_1,
  \sigma_2 \in \C0 \left( \rpic; \mathcal{X}_\alpha \right)$,
  by~\textbf{(V1)}:
  \begin{eqnarray*}
    \norma{\nabla_x^2 V(\sigma_1)}_{\L\infty(I;\L1(\reali^N;\reali^{N\times N \times N}))}
    & \leq &  
    C \,,
    \\
    \norma{V(\sigma_1)-V(\sigma_2)}_{\L\infty(I\times\reali^N;\reali )}
    & \leq &
    C \, \norma{\sigma_1-\sigma_2}_{\L\infty(I;\L1(\reali^N;\reali))}\,,
    \\
    \norma{\div\left(V(\sigma_1) - V(\sigma_2)\right)}_{\L\infty(I;\L1(\reali^N;\reali))}
    & \leq &
    C \, \norma{\sigma_1-\sigma_2}_{\L\infty(I;\L1(\reali^N;\reali))}\,.
  \end{eqnarray*}
  Applying~\cite[Theorem~2.6]{ColomboMercierRosini}, we get
  \begin{eqnarray*}
    \norma{(\rho_1-\rho_2)(t)}_{\L1}
    & \leq &
    C t e^{C't}\tv(\rho_o) \norma{\sigma_1-\sigma_2}_{\L\infty([0,t]; \L1)}
    \\
    & &
    \quad +  
    C^2NW_N \int_0^t  (t-s)e^{C'(t-s)} \,  \d{s}
    \, \norma{\sigma_1-\sigma_2}_{\L\infty([0,t];\L1)} 
    \\
    & &
    \quad +
    \int_0^t C \, e^{ C'(t-s)} \,
    \norma{(\sigma_1-\sigma_2)(s)}_{\L1} \,  \d{s}\, .
  \end{eqnarray*}
  So that
  \begin{displaymath}
    \norma{\mathcal{Q} (\sigma_1) - \mathcal{Q}(\sigma_2)}_{\L\infty(I;\L1)}
   \leq
    C T e^{C'T}
    \left[ \tv(\rho_o) + N W_N C T + 1\right] 
    \norma{\sigma_1-\sigma_2}_{\L\infty(I; \L1)} \, .
  \end{displaymath}
  Here we introduce the function $f(T) = C T e^{C'T} \left[
    \tv(\rho_o) + N W_N C T + 1\right]$ and we remark that $f(T) \to
  0$ when $T \to 0$. Choose now $T_1 > 0$ so that $f(T_1) =
  \frac{1}{2}$.  Banach Contraction Principle now ensures the
  existence and uniqueness of a solution to~(\ref{eq:SCL}) on
  $[0,T_1]$ in the sense of Definition~\ref{def:sol}.

  Iterate this procedure up to the interval $[T_{n-1}, T_n]$ and
  obtain the contraction estimate
  \begin{eqnarray*}
    & &
    \norma{\mathcal{Q} (\sigma_1) - \mathcal{Q}(\sigma_2)}_{\L\infty([T_{n},T_{n+1}];\L1)}
    \\
    & \leq &
    C(T_{n+1}-T_n)e^{C'(T_{n+1}-T_n)}
    \left[ 
      \tv \left(\rho(T_n) \right) + N W_N C ( T_{n+1} - T_n) +1
    \right]
    \\
    & &
    \times 
    \norma{\sigma_1-\sigma_2}_{\L\infty([T_{n},T_{n+1}]; \L1)}
    \\
    & \leq &
    \left[ 
      \tv(\rho_o)\,e^{C'T_n}+C'T_n e^{C'T_n} +NW_N C(T_{n+1}-T_n)
      +1\right]
    \\
    & &
    \times   C'(T_{n+1}-T_n)e^{C'(T_{n+1}-T_n)}
    \norma{\sigma_1-\sigma_2}_{\L\infty([T_{n},T_{n+1}]; \L1)}
  \end{eqnarray*}
  where we used the bounds on $\tv \left( \rho(T_n) \right)$ and
  $\norma{\rho(T_n)}_{\L\infty(\reali^N;\reali)}$ provided by
  Lemma~\ref{lem:Infty} associated to the conditions~\textbf{(A)}
  and~\textbf{(V1)}. We may thus extend the solution up to time
  $T_{n+1}$, where we take
  \begin{eqnarray*}
    \left[\tv(\rho_o)\,e^{C'T_n}+CT_n e^{C'T_n} +NW_N C(T_{n+1}-T_n)
      +1\right]
    C(T_{n+1}-T_n)e^{C'(T_{n+1}-T_n)} 
    & = &
    \frac{1}{2}\,.
  \end{eqnarray*}
  If the sequence $(T_n)$ is bounded, then the left hand side above
  tends to $0$, whereas the right hand side is taken equal to $1/2 >
  0$. Hence, the sequence $(T_n)$ is unbounded and the solution
  to~(\ref{eq:SCL}) is defined on all $\rpic$.

  {(S\ref{it:s4})} follows from Lemma~\ref{lem:Infty} associated to
  the assumption~\textbf{(V1)} on $V$ that allows to satisfy the
  hypotheses on $w$.

  {(S\ref{it:s3})} is obtained in the same way as
  {(S\ref{it:s1})}. Note that the Lipschitz constant obtained by such
  a way is depending on time.

  The bound~{(S\ref{it:s2})} follows from Lemma~\ref{lem:Infty}, point
  2, that gives us
  \begin{displaymath}
    \tv\left( \rho(t) \right) 
    \leq 
    \tv(\rho_o) e^{C't} + N W_N Cte^{C't} \norma{\rho_o}_{\L\infty} \,.
  \end{displaymath}

  \smallskip

 \noindent\textbf{(B)}: \quad Repeat the proof of
 Theorem~\ref{thm:main} and, with the notation therein, note that if
 we find a sequence $(\alpha_n)$ such that $\sum_n
 T(\alpha_n,\alpha_{n+1})=+\infty$ where $T(\alpha,\beta) =
 \left[\ln\left(\beta/\alpha\right)\right] / C(\beta)$, then the
 solution is defined on the all $\rpic$. It is immediate to check
 that~\textbf{(B)} implies that
 \begin{displaymath}
   \sum_{n=1}^k T(\alpha_n,\alpha_{n+1})
   \geq
   \left( \norma{C}_{\L\infty(\rpic; \rpic)} \right)^{-1}
   \ln \alpha_k
   \to +\infty
   \quad
   \mbox{as } k \to +\infty
 \end{displaymath}
 completing the proof.
\end{proofof}

\begin{proofof}{Proposition~\ref{prop:prop}}
  The bounds of $\rho$ in $\W1\infty$ and $\W11$ follow
  from~\ref{it:Infty:4} in Lemma~\ref{lem:Infty}, the hypotheses being
  satisfied thanks to~\textbf{(V2)}. The bound in $\W21$ comes
  from~\ref{it:Infty:5} in Lemma~\ref{lem:Infty}, the hypotheses being
  satisfied thanks to~\textbf{(V3)}.
\end{proofof}

\subsection{Weak G\^ateaux Differentiability}

First of all, if $r_o\in (\L\infty \cap \L1)(\reali^N;\reali)$ and
$\rho \in \L\infty \bigl(I_{\mathrm{ex}}; (\W11 \cap \W1\infty)
(\reali^N;\reali) \bigr)$, we prove that the
equation~(\ref{eq:linear}) admits a unique solution $r\in \L\infty
\bigl( I_{\mathrm{ex}}; \L1(\reali^N;\reali) \bigr)$ continuous from the
right.

\begin{proofof}{Proposition~\ref{prop:existence_r}}
  We use here once again Lemma~\ref{lem:Both} in order to get an
  expression of the Kru\v{z}kov solution for~(\ref{eq:simple}).

  We assume now that $\rho\in \C0 \bigl(I_{\mathrm{ex}}; (\W1\infty \cap
  \W11)(\reali^N; \reali) \bigr)$ and we define $w = V(\rho)$; we also
  set, for all $s \in \L\infty \bigl( I_{\mathrm{ex}};
  \L1(\reali^N;\reali) \bigr)$, $R = \div \left(\rho DV(\rho)(s)
  \right)$. Thanks to the assumptions on $\rho$ and~\textbf{(V4)}, we
  obtain $R\in \L\infty \bigl( I_{\mathrm{ex}}; \L1(\reali^N;\reali)
  \bigr) \cap \L\infty(I_{\mathrm{ex}} \times \reali^N; \reali)$. Let
  $\epsilon \in \pint{I}_{\mathrm{ex}}$. Then, on
  $[0,T_{\mathrm{ex}}-\epsilon]$ we can apply Lemma~\ref{lem:Both}
  giving the existence of a Kru\v zkov solution to
  \begin{displaymath}
    \pt_t r+\div(rw)=R \,, \quad 
    r(x,0) = r_{o} \in  (\L\infty \cap \L1)(\reali^N;\reali) \,.
  \end{displaymath}
  Let $T \in [0,T_{\mathrm{ex}}-\epsilon]$ and $I = \left[0 , T
  \right[$. We denote $Q$ the application that associates to $s\in
  \L\infty \bigl(I;\L1(\reali^N;\reali) \bigr)$ continuous from the
  right in time, the Kru\v{z}kov solution $r \in \L\infty
  \bigl(I;\Lloc1(\reali^N;\reali) \bigr)$ continuous from the right in
  time of~(\ref{eq:simple}) with initial condition $r_o\in (\L\infty
  \cap \L1)(\reali^N;\reali)$, given by Lemma~\ref{lem:Both}. That is
  to say
  \begin{eqnarray*}
    Q 
    &\colon& 
    s\mapsto r(t,x)  = r_o\left( X(0; t,x) \right)
    \exp
    \left(
      -\int_0^t \div V(\rho) \left( \tau, X(\tau; t,x) \right) \d{\tau}
    \right)
    \\
    & & 
    -\int_0^t \div \left(\rho DV(\rho) (s) \right) 
    \left(\tau, X(\tau;t,x) \right)
    \exp \left(
      -\int_\tau^t \div V(\rho)\left(u, X (u;t,x) \right) \d{u}
    \right) \,\d{\tau} .
  \end{eqnarray*}
  Let us give some bounds on $r$. The representation of the
  solution~(\ref{eq:r}) allows indeed to derive a $\L\infty$ bound on
  $r$. For all $t\in I$, thanks to~\textbf{(V1)} and~\textbf{(V4)} we
  get, with $C = C \left(
    \norma{\rho}_{\L\infty([0,T_{\mathrm{ex}}-\epsilon]\times\reali^N;\reali)}
  \right)$,
  \begin{displaymath}
    \norma{r(t)}_{\L\infty} \leq\norma{r_o}_{\L\infty}e^{Ct}
    +
    t e^{Ct} 
    \norma{\rho}_{\L\infty([0,t];\W1\infty)}
    \norma{DV(\rho)}_{\W1\infty}
    \norma{s}_{\L\infty([0,t],\L1)}\,.
  \end{displaymath}
  The same expression allows also to derive a $\L1 $ bound on $r(t)$
  \begin{displaymath}
    \norma{r(t)}_{\L1} 
    \leq
    \norma{r_o}_{\L1}e^{Ct}
    +
    t e^{Ct}\norma{\rho}_{\L\infty([0,t];\W11)}
    \norma{DV(\rho)}_{\W1\infty} \norma{s}_{\L\infty([0,t],\L1)}\,.
  \end{displaymath}

  Now, we want to show that $Q$ is a contraction.  We use once again
  the assumption~\textbf{(V4)}. For all $s_1, s_2 \in \L\infty \bigl(
  I; (\L1 \cap \BV)(\reali^N;\reali) \bigr)$ continuous from the
  right, we have
  \begin{displaymath}
    \norma{\div\left( \rho
        DV(\rho)(s_1-s_2)\right)}_{\L1(\reali^N;\reali)} 
    \leq 
    C
    \norma{\rho}_{\W11(\reali^N;\reali)}\norma{s_1-s_2}_{\L1(\reali^N;\reali)}\, .
  \end{displaymath}
  Thus, we get:
  \begin{eqnarray*}
    & &
    \norma{Q(s_1)-Q(s_2)}_{\L\infty(I;\L1)}
    \\
    &\leq &
    C
    \norma{\rho}_{\L\infty(I;\W11)} 
    \norma{s_1-s_2}_{\L\infty(I;\L1)}  
    \int_0^T \exp 
    \left(
      (T-\tau)\norma{\div V(\rho)}_{\L\infty}    
    \right) \d{\tau}
    \\
    & \leq & 
    (e^{CT}-1)
    \norma{\rho}_{\L\infty([0,T_{\mathrm{ex}}-\epsilon];\W11)} 
    \norma{s_1-s_2}_{\L\infty(I;\L1)}\,.
  \end{eqnarray*}
  Then, for $T$ small enough, can apply the Fixed Point Theorem, that
  gives us the existence of a unique Kru\v{z}kov solution to the
  problem. Furthermore, as the time of existence does not depend on
  the initial condition, we can iterate this procedure to obtain
  existence on the interval $[0,T_{\mathrm{ex}}-\epsilon]$. Finally, as
  this is true for all $\epsilon \in \pint{I}_{\mathrm{ex}}$, we obtain
  the same result on the all interval $I_{\mathrm{ex}}$.

  The $\L1$ bound follows from~(\ref{eq:r}). Let $T \in I_{\mathrm{ex}}$
  and $t\in I$, then for a suitable $C = C\left(
    \norma{\rho}_{\L\infty(I\times \reali^N;\reali)} \right)$
  \begin{displaymath}
    \norma{r(t)}_{\L1}
    \leq 
    \norma{r_o}_{\L1} e^{Ct} 
    + 
    \norma{\rho}_{\L\infty(I;\W11)} \,
    \norma{\div DV(\rho)}_{\L\infty}  \,
    \int_0^t \norma{r(\tau)}_{\L1} \d{\tau} \, .
  \end{displaymath}
  A use of~\textbf{(V4)} and an application of Gronwall Lemma gives
  \begin{eqnarray*}
    \norma{r(t)}_{\L1}
    &\leq & 
    e^{Ct}
    e^{K\norma{\rho}_{\L\infty(I,\W11)}t} 
    \, \norma{r_o}_{\L1} \,,
  \end{eqnarray*}
  where
  $K=K\left(\norma{\rho}_{\L\infty(I\times\reali^N;\reali)}\right)$ is
  as in~\textbf{(V4)}.

  The $\L\infty $ bound comes from the same representation
  formula. Indeed, for $T\in I_{\mathrm{ex}}$ and $ t \in I$ we have
  \begin{eqnarray*}
    \norma{r(t)}_{\L\infty}
    \leq e^{Ct} 
    \norma{r_o}_{\L\infty}
    + 
    \norma{\rho}_{\L\infty(I;\W1\infty)} \,
    \norma{\div DV(\rho)}_{\L\infty}  \,
    \int_0^t \norma{r(\tau)}_{\L1} \d{\tau} \, .
  \end{eqnarray*}
  Then, the last $\norma{r(\tau )}_{\L1}$ is bounded just as above. We
  get
  \begin{displaymath}
    \norma{r(t)}_{\L\infty}
    \leq
    e^{Ct} \norma{r_o}_{\L\infty} 
    +
    K t e^{2Ct}
    e^{K\norma{\rho}_{\L\infty(I,\W11)}t}
    \norma{r_o}_{\L1}
    \norma{\rho}_{\L\infty(I,\W1\infty)} \,.
  \end{displaymath}

  Finally, we get a $\W11$ bound using the expression of the solution
  given by Lemma~\ref{lem:Both}. Indeed, assuming in
  addition~\textbf{(V2)} and~\textbf{(V4)}, we get
  \begin{eqnarray*}
    \norma{\nabla r(t)}_{\L1(\reali^N;\reali)}
    & \leq & 
    e^{2Ct} \norma{\nabla r_o}_{\L1}
    +
    C t e^{2Ct} \norma{r_o}_{\L1}
    \\
    & & 
    + 
    K (1 +  Ct) e^{2Ct}
    \norma{\rho}_{\L\infty(I;\W21)} \int_0^t \norma{r(\tau)}_{\L1} \d{\tau}
    \\
    & \leq &
    e^{2Ct} \norma{\nabla r_o}_{\L1}
    +
    C t e^{2Ct} \norma{r_o}_{\L1}
    \\
    & & 
    + 
    K t (1 +  Ct) e^{3Ct} 
    e^{K\norma{\rho}_{\L\infty(I,\W11)}t} 
    \norma{r_o}_{\L1}
    \norma{\rho}_{\L\infty(I;\W21)}
    \,.
  \end{eqnarray*}
  Hence, denoting $C' = \max \{C, K\norma{\rho}_{\L\infty(I,\W11)}\}$,
  we obtain
  \begin{displaymath}
    \norma{r(t)}_{\W11}
    \leq 
    \norma{r_o}_{\W11}(1+C't) e^{2C't}
    + 
    K t (1 +  Ct) e^{4C't} 
    \norma{r_o}_{\L1}
    \norma{\rho}_{\L\infty(I;\W21)}
    \,.
  \end{displaymath}

  The full continuity in time follows
  from~\cite[Remark~2.4]{ColomboMercierRosini} and~\textbf{(V1)},
  \textbf{(V4)}, since $r(t)\in \W11(\reali^N;\reali)$ implies that
  $r(t)\in \BV(\reali^N;\reali)$ with $\tv(r(t))=\norma{\nabla_x
    r(t)}_{\L1(\reali^N;\reali)}$.
\end{proofof}

Now, we can address the question of weak G\^ateaux differentiability
of the semigroup giving the solution to~(\ref{eq:SCL}).

\begin{proofof}{Theorem~\ref{thm:Diff}} Let $\alpha, \beta >0$ with
  $\beta>\alpha$ and $h\in [0,h^*]$ with $h^*$ small enough so that
  $\beta>\alpha(1+h^*)$.  Fix $\rho_o, r_o \in \mathcal{X}_\alpha$.
  Thanks to Theorem~\ref{thm:main}, we get the weak entropy solution
  $\rho \in \C0([0,T(\alpha,\beta)]; \mathcal{X}_{\beta})$
  of~(\ref{eq:SCL}) with initial condition $\rho_o$ and $\rho_h\in
  \C0([0,T(\alpha(1+h),\beta)];\mathcal{X}_{\beta})$ of~(\ref{eq:SCL})
  with initial condition $\rho_o+hr_o$. Note that
  \begin{displaymath}
    T(\alpha(1+h),\beta)
    =
    \frac{\ln(\beta/(\alpha(1+h)))}{C(\beta)}
    =
    T(\alpha,\beta)-\frac{\ln(1+h)}{C(\beta)}
    \leq
    T(\alpha,\beta)
  \end{displaymath}
  and $T(\alpha(1+h),\beta)$ goes to $T(\alpha,\beta)$ when $h$ goes
  to 0. In particular, both solutions are defined on the interval
  $[0,T(\alpha(1+h^*),\beta)]$.

  By Theorem~\ref{thm:main}, point 2, the sequence $\left(\frac{\rho_h
      - \rho}{h}(t)\right)_{h\in [0,h^*]}$ is bounded in $\L1$ for all
  $t\in [0, T(\alpha(1+h^*),\beta)]$. By Dunford--Pettis Theorem, it
  has a weakly convergent subsequence, see~\cite{Brezis}. Thus, there
  exists $r\in\L1(\reali^N;\reali)$ such that
  \begin{displaymath}
    \displaystyle \frac{\rho_h- \rho}{h}(t)\rightharpoonup_{h\to 0} r(t)
    \textrm{ weakly in } \L1.
  \end{displaymath}
  Write now the definition of weak solution for $\rho$, $\rho_h$. Let
  $\phi \in \Cc\infty ([0,T(\alpha(1+h^*),\beta)] \times \reali^N;
  \reali)$
  \begin{eqnarray*}
    \int_\rpis \int_{\reali^N}
    \left(
      \rho \pt_t \phi
      +
      \left( \rho V(\rho) \right) \cdot \nabla_x \phi
    \right)\dx\, \dt
    & = &
    0\,;
    \\
    \int_\rpis \int_{\reali^N}
    \left(
      \rho_h \pt_t \phi
      +
      \left( \rho_h V(\rho_h) \right)
      \cdot \nabla_x \phi
    \right)
    \dx\, \dt
    & = &
    0\,.
  \end{eqnarray*}
  Now, use~\textbf{(V4)} and write, for a suitable function $\epsilon
  = \epsilon(\rho,\rho_h)$,
  \begin{displaymath}
    V(\rho_h)
    =
    V(\rho)
    +
    \mathrm{D}V(\rho) (\rho_h - \rho)
    +
    \epsilon(\rho, \rho_h)\,,
  \end{displaymath}
  with $\norma{\epsilon(\rho,\rho_h)}_{\L\infty(\reali^N;\reali)}\leq
  K(2\beta) \left(\norma{ \rho_h - \rho}_{\L1(\reali^N;\reali)}
  \right)^2$.  Then,
  \begin{displaymath}
    \rho V(\rho) -  \rho_h V(\rho_h)
    = 
    (\rho-\rho_h) V(\rho)
    +
    \rho \mathrm{D}V(\rho)(\rho-\rho_h)
    +
    (\rho - \rho_h) \mathrm{D}V(\rho) (\rho-\rho_h)
    -
    \rho_h \epsilon(\rho,\rho_h).
  \end{displaymath}
  Consequently,
  \begin{eqnarray*}
    \int_{\reali_{+}^*} \int_{\reali^N}
    \left[
      \frac{\rho- \rho_h}{h} 
      \pt_t \phi
      +
      \left(
        \frac{\rho- \rho_h}{h} V(\rho)
        +
        \rho \mathrm{D}V(\rho)
        \left(\frac{\rho- \rho_h}{h}\right)
      \right.
    \right.
    \\
    \left.\left. +\frac{\rho- \rho_h}{h}
        \mathrm{D}V(\rho) (\rho-\rho_h)-\rho_h \frac{\epsilon (
          \rho,\rho_h)}{h}\right)\cdot \nabla_x \phi \right]\dx\, \dt
    & = &
    0\,.
  \end{eqnarray*}
  Using~\textbf{(V4)}, $\rho(t) \in \mathcal{X}_{\beta}$ and the
  estimate on $\epsilon$ we obtain for all $t \in [0, T(
  \alpha(1+h^*), \beta) ]$:
  \begin{eqnarray*}
    & &
    \int_{\reali^N} 
    \modulo{\nabla_{x}\varphi }
    \modulo{\frac{\rho - \rho_h}{h} \mathrm{D}V(\rho)
      (\rho-\rho_h)-\rho_h \frac{\epsilon ( \rho, \rho_h)}{h}} \dx 
    \\
    & \leq &
    K(2\beta) 
    \int_{\reali^N}  \left(
      \frac{\modulo{\rho-\rho_h}}{h}
      +
      \beta \frac{\norma{\rho-\rho_h}_{\L1}}{h}
    \right)
    \norma{\rho- \rho_h}_{\L1}\, 
    \modulo{\nabla_{x}\varphi } \dx,
  \end{eqnarray*}
  and since $\frac{\rho_h-\rho}{h}$ is bounded in $\L\infty \left(
    [0,T(\alpha(1+h^*),\beta)]; \L1(\reali^N;\reali) \right)$ and
  $\frac{\rho_h-\rho}{h} (t) \wto_{h\to 0} r(t)$ in $\L1$, then we can
  apply the Dominated Convergence Theorem. We get:
  \begin{eqnarray*}
    \int_\reali \int_{\reali^N}
    & &
    \left[
      r \pt_t  \phi
      +
      \left(
        r V(\rho)+\rho \mathrm{D} V(\rho)
        (r)
      \right)
      \cdot \nabla_x \phi
    \right]
    \dx\, \dt =0\, .
  \end{eqnarray*}
  That is to say that $r$ is a weak solution to~(\ref{eq:linear}) with
  initial condition $r_o$. As this is true for all $h^*$ small enough,
  finally we obtain a solution on the all interval $[0,T(\alpha,
  \beta)[$.  Hence we conclude that $\rho\in \C0(I_{\mathrm{ex}} \times
  \reali^N;\reali)$ implies $r$ defined on $I_{\mathrm{ex}}$.
\end{proofof}

In the proof just above, we can not conclude to the uniqueness of the
weak G\^ateaux derivative as we do not know if the weak solution is
unique In particular, we don't know if the derivative is continuous.

We assume now that the assumptions~\textbf{(V4)} and~\textbf{(V5)} are
satisfied by $V$. We want to show that with these hypotheses, we have
now strong convergence in $\L1$ to the Kru\v{z}kov solution
of~(\ref{eq:linear})

\begin{proofof}{Theorem~\ref{thm:Diff2}}
  Let $\alpha, \beta >0$ with $\beta>\alpha$, and $h\in [0,h^*]$ with
  $h^*$ small enough so that $\beta>\alpha(1+h^*)$.  Let us denote
  $T(h)=T(\alpha(1+h),\beta)$ for $h\in [0,h^*]$ the time of existence
  of the solution of~(\ref{eq:SCL}) given by Theorem~\ref{thm:main}.

  Fix $\rho_o \in ( \W1\infty \cap \W21) (\reali^N;[0,\alpha])$, $ r_o
  \in (\L\infty \cap \W11) (\reali^N; [0,\alpha])$. Let $\rho$,
  respectively $ \rho_h$, be the weak entropy solutions
  of~(\ref{eq:SCL}) given by Theorem~\ref{thm:main} with initial
  condition $\rho_o$, respectively $\rho_o+h r_o$. Note that these
  both solutions are in $\C0\bigl( [0,T(h^*)]; \L1(\reali^N;\reali)
  \bigr)$.  Furthermore, under these hypotheses for $\rho_o$ and
  $r_o$, we get thanks to Proposition~\ref{prop:prop} that the
  corresponding solutions $\rho$ and $\rho_h$ of~(\ref{eq:SCL}) are in
  $\C0 \bigl( [0,T(h^*)]; (\W1\infty \cap \W21) (\reali^N ; [0,\beta])
  \bigr)$, condition~\textbf{(V3)} being satisfied. Hence, we can now
  introduce the Kru\v{z}koz solution $r \in \C0 \bigl( [0,T(h^*)[;
  \L1(\reali^N;\reali) \bigr)$ of~(\ref{eq:linear}), whose existence
  is given in this case by Proposition~\ref{prop:existence_r}. Note
  that, $r_o$ being in $\W11(\reali^N; \reali)$ and $\rho \in \L\infty
  \left( [0,T(h^*)] ; \W21(\reali^N; \reali) \right)$
  and~\textbf{(V2), (V4)} being satisfied, $r(t)$ is also in
  $\W11(\reali^N; \reali)$ for all $t \in \left[ 0, T(h^*)\right[$
  thanks to the $\W11$ bound of Proposition~\ref{prop:existence_r}.

  % First, by~\textbf{(S3)}, the sequence $\frac{\rho_h - \rho}{h}$ is
  % bounded in $\L1$. By Dunford--Pettis Theorem, it has a weakly
  % convergent subsequence, see~\cite{Brezis}.

  Let us denote $z_h=\rho+hr$. We would like to compare $\rho_h$ and
  $z_h$ thanks to~\cite[Theorem 2.6]{ColomboMercierRosini}.  A
  straightforward computation shows that $z_h$ is the solution to the
  following problem,
  \begin{displaymath}
    \left\{ \begin{array}{l}
        \pt_t z_h + \div\left(z_h \left(V(\rho)+hDV(\rho)(r) \right) \right)
        =
        h^2 \div \left( r \, DV(\rho)(r) \right)\, , 
        \\
        z_h(0) = \rho_o + h r_o \in \mathcal{X}_{\alpha(1+h)} \, .
      \end{array}
    \right.
  \end{displaymath}
  Note that the source term being in $\C0 \bigl( \left[0, T(h^*)
  \right[; \L1(\reali^N;\reali) \bigr)$, and the flow being regular,
  we can apply to this equation Lemma~\ref{lem:Both} that gives
  existence of a Kru\v zkov solution.

  As in the proof of Lemma~\ref{lem:Both}, we make here the remark
  that~\cite[Theorem 2.6]{ColomboMercierRosini} can be used with the
  second source term in $\C0 \left( \left[0, T(h^*) \right[;
    \L1(\reali^N; \reali) \right)$ and the flow $\C2$ in space and
  only $\C0$ in time. Besides, we also use the same \emph{slight
    improvement} as in the proof of Lemma~\ref{lem:Infty}, taking the
  $\L\infty$ norm in the integral term only in space, keeping the time
  fixed.  We get, with $\kappa_o = N W_N (2N+1) \norma{\nabla_x
    V(\rho_h)}_{\L\infty([0,T(h^*)]\times\reali^N;\reali)}$ and
  $\kappa = 2 N\norma{\nabla_x
    V(\rho_h)}_{\L\infty([0,T(h^*)]\times\reali^N;\reali)}$, for some
  $T \in [0,T(h^*)]$,
  \begin{eqnarray*}
    & & 
    \norma{\rho_h-z_h}_{\L\infty(I;\L1)}
    \\
    & \leq & 
    T e^{\kappa_o T}\tv (\rho_o +hr_o) 
    \norma{V(\rho_h)-V(\rho)-h DV(\rho)(r)}_{\L\infty([0,T(h^*)]\times\reali^N;\reali^N)}
    \\
    & &
    + 
    N W_N \int_0^T (T-t) e^{\kappa_o(T-t)} \int_{\reali^N}
    \norma{\rho_h(t)}_{\L\infty} 
    \norma{\nabla_x \div V(\rho_h)} \dx \dt
    \\
    & &
    \qquad \qquad
    \times 
    \norma{V(\rho_h)-V(\rho)-h DV(\rho)(r)}_{\L\infty([0,T(h^*)]\times\reali^N;\reali^N)}
    \\
    & & 
    + h^2\int_0^T e^{\kappa(T-t)} \int_{\reali^N} 
    \modulo{\div \left(rDV(\rho)(r) \right)} \dx \, \dt
    \\
    & &
    + \int_0^T e^{\kappa(T-t)} \int_{\reali^N} 
    \modulo{\div \left(V(\rho_h)-V(\rho)-h DV(\rho)(r) \right)} \dx \,
    \dt\\&&\qquad\qquad\times
    \max_{t\in[0,T]} \left\{ 
      \norma{\rho_h (t)}_{\L\infty}, \norma{z_h(t)}_{\L\infty} 
    \right\} .
  \end{eqnarray*}
  Then, setting $C = C(\beta)$ and $K=K(2\beta)$, we use:
  \begin{itemize}
  \item the bound of $\rho$ and $\rho_h$ in $\L\infty$ given by
    Lemma~\ref{lem:Infty}
    \begin{eqnarray*}
      \norma{\rho(t)}_{\L\infty}
      \leq 
      \norma{\rho_o}_{\L\infty} e^{Ct}
      \leq 
      \beta
      \quad \mbox{ and } \quad
      \norma{\rho_h(t)}_{\L\infty}
      \leq
      \norma{\rho_o+h r_o}_{\L\infty} e^{Ct}
      \leq
      \beta\,;
    \end{eqnarray*}
  \item the properties of $V$ given in~\textbf{(V1)} to get
    \begin{displaymath}
      \norma{\nabla_x \div V(\rho_h)}_{\L\infty(\rpic\times\reali^N;\reali)} 
      \leq C 
      \quad \mbox{ and } \quad
      \norma{\nabla_x \div V(\rho_h)}_{\L\infty(I;\L1(\reali^N;\reali))} 
      \leq C \,;
    \end{displaymath}
  \item the property~\textbf{(V4)}, respectively~\textbf{(V5)}, to get
    \begin{eqnarray*}
      & & 
      \norma{V(\rho_h)-V(\rho)-h
        DV(\rho)(r)}_{\L\infty(I\times\reali^N;\reali)}
      \\
      & \leq &
      K \, \left(
        \norma{\rho_h- \rho}_{\L\infty(I;\L1(\reali^N;\reali))}^2 
        +    
        \norma{\rho_h-z_h}_{\L\infty(I;\L1(\reali^N;\reali))}\right)
      \,, \quad \mbox{respectively}
      \\
      & &
      \norma{\div\left(V(\rho_h)-V(\rho)-h
          DV(\rho)(r)\right)}_{\L\infty(I;\L1(\reali^N;\reali))}
      \\
      &\leq &
      K \, \left(
        \norma{\rho_h- \rho}_{\L\infty(I;\L1(\reali^N;\reali))}^2  
        + 
        \norma{\rho_h-z_h}_{\L\infty(I;\L1(\reali^N;\reali))}\right)\, ;
    \end{eqnarray*}
  \item the property~\textbf{(V4)} to get
    \begin{displaymath}
      \norma{\div\left( r DV(\rho)(r)\right)}_{\L1} 
      \leq 
      K \,      \norma{r}_{\W11} 
      \norma{r}_{\L1} \, .
    \end{displaymath}
  \end{itemize}
  Gathering all these estimates, denoting $C'=NW_N(2N+1)C$, we obtain
  \begin{eqnarray*}
    & &    \norma{\rho_h-z_h}_{\L\infty(I;\L1)}\\
    & \leq & 
    T e^{C'T}
    \left( 
      \tv(\rho_o+hr_o) 
      +
      N W_N CT \beta
    \right) 
    K \left(
      \norma{\rho_h- \rho}_{\L\infty(I;\L1)}^2 
      +
      \norma{\rho_h-z_h}_{\L\infty(I;\L1)} 
    \right)
    \\
    & &
    + 
    h^2 K T e^{C' T} \norma{r}_{\L\infty(I;\W11)}
    \norma{r}_{\L\infty(I;\L1)} 
    \\
    & &
    + 
    \left(
      \beta
      + 
      h \sup_{t\in I} \norma{r(t)}_{\L\infty}\right) T e^{C'T}
    K\left(
      \norma{\rho_h- \rho}_{\L\infty(I;\L1)}^2 
      + 
      \norma{\rho_h-z_h}_{\L\infty(I;\L1)} 
    \right) \,.
  \end{eqnarray*}
  Then, dividing by $h$ and introducing
  \begin{eqnarray*}
    F_h(T)
    & = &
    K Te^{C'T} 
    \left[
      \tv(\rho_o)+h\tv(r_o) 
      +
      N W_N C T \beta 
      +
      \beta
      +
      h\norma{r(t)}_{\L\infty}
    \right] \,,
  \end{eqnarray*}
  we obtain
  \begin{eqnarray*}
    \norma{\frac{\rho_h-z_h}{h}}_{\L\infty(I;\L1)}
    & \leq &
    F_h(T)
    \left[
      \norma{\rho_h-\rho}_{\L\infty(I;\L1)}
      \norma{\frac{\rho_h - \rho}{h}}_{\L\infty(I;\L1)}
      + 
      \norma{\frac{\rho_h-z_h}{h}}_{\L\infty(I;\L1)} 
    \right]
    \\ 
    & &
    +
    h K T e^{C' T }
    \norma{r}_{\L\infty(I;\W11)} 
    \norma{r}_{\L\infty(I;\L1)}  \, .
  \end{eqnarray*}
  Note that $F_h$ is a function that vanishes in $T=0$ and that
  depends also on $\rho_o$, $r_o$ and $h$.  Hence, we can find $\bar
  T\leq T(h^*)$ small enough such that $ F_{h^*}(\bar T) \leq
  1/2$. Furthermore, $F_h(T)$ is increasing in $h$ consequently,
  $h\leq h^*$ implies $F_h(T)\leq F_{h^*}(T)$.  Noticing moreover that
  $\norma{\frac{\rho_h - \rho}{h}}_{\L\infty(I;\L1)}$ has a uniform
  bound $M$ in $h$ by 2.~in Theorem~\ref{thm:main}, we get for $T\leq
  \bar T$
  \begin{eqnarray*}
    \frac{1}{2}\norma{\frac{\rho_h-\rho}{h}-r}_{\L\infty(I;\L1)}
    & = & 
    \frac{1}{2}\norma{\frac{\rho_h-z_h}{h}}_{\L\infty(I;\L1)}\\
    &\leq &
    \frac{M}{2}  \norma{\rho_h-\rho}_{\L\infty(I;\L1)}
%     \\
%     & & 
    + 
    h KT e^{C'T} 
    \norma{r}_{\L\infty(I,\W11)} 
    \norma{r}_{\L\infty(I;\L1)}\, .
  \end{eqnarray*}
  The right side above goes to 0 when $h\to 0$, so we have proved the
  G\^ateaux differentiability of the semigroup $S$ for small
  time. Finally, we iterate like in the proof of
  Theorem~\ref{thm:main} in order to have existence on the all
  interval $[0,T(h^*)]$. Let $T_1$ be such that $F_{h^*}(T_1)=1/2$ and
  assume $T_1<T(h^*)$. If we assume the G\^ateaux differentiability is
  proved until time $T_n\leq T(h^*)$, we make the same estimate on
  $[T_n,T_{n+1}]$, $T_{n+1}$ being to determine. We get
  \begin{eqnarray*}
    & &
    \norma{\rho_h-z_h}_{\L\infty([T_n,T_{n+1}];\L1)}
    \\
    & \leq & 
    (T_{n+1}-T_n) e^{C' (T_{n+1}-T_n)}
    \left( 
      \tv(\rho_h(T_n)) 
      +
      N W_N C (T_{n+1}-T_n)
      \beta
    \right) 
    \\
    & &
    \qquad\times
    K\left(
      \norma{\rho_h- \rho}_{\L\infty([ T_{n},T_{n+1}];\L1)}^2 
      +
      \norma{\rho_h-z_h}_{\L\infty([T_{n},T_{n+1}];\L1)} 
    \right)
    \\
    & &
    + 
    h^2 K  (T_{n+1}-T_n) e^{C'(T_{n+1}-T_n)} 
    \norma{r}_{\L\infty([T_{n},T_{n+1}];\W11)}
    \norma{r}_{\L\infty([T_{n},T_{n+1}];\L1)}
    \\
    & &
    + 
    \left(\beta+h\sup_{[T_n,T_{n+1}]}\norma{r(t)}_{\L\infty}\right)  (T_{n+1}-T_n) e^{C' (T_{n+1}-T_n)}
    \\
    & &
    \qquad \times 
    K\left(
      \norma{\rho_h- \rho}_{\L\infty([T_{n},T_{n+1}];\L1)}^2 
      + 
      \norma{\rho_h-z_h}_{\L\infty([T_{n},T_{n+1}];\L1)} 
    \right) \,.
  \end{eqnarray*}
  Then, we divide by $h$ and we introduce, for $T\geq T_n$
  \begin{eqnarray*}
    F_{h,n}(T)
    & = &
    K (T-T_n)e^{C'(T-T_n)} 
    \left[
      (\tv(\rho_o)+h\tv(r_o))e^{CT_n}
      +
      \beta C'T_n e^{C'T_n} 
    \right. 
    \\
    & &
    \left. 
      +
      N W_N C(T-T_n) \beta
      +
      \beta + h \sup_{[T_n,T_{n+1}]} \norma{r(t)}_{\L\infty}
    \right]\,.
  \end{eqnarray*}
  We define $T_{n+1}>T_n$ such that
  $F_{h,n}(T_{n+1})=\frac{1}{2}$. This is possible since $F_{h,n}$
  vanishes in $T=T_n$ and increases to infinity when $T\to\infty$.
  Hence, as long as $T_{n+1} \leq T(h^*)$, we get
  \begin{eqnarray*}
    & &
    \norma{\frac{\rho_h-\rho}{h}-r}_{\L\infty([T_n,T_{n+1}];\L1)}
    \\
    &\leq &
    K M \norma{\rho_h-\rho}_{\L\infty([T_n,T_{n+1}];\L1)}
    \\
    & & 
    + 
    2 h K(T_{n+1}-T_n) e^{C'(T_{n+1}-T_n)} 
    \norma{r}_{\L\infty([T_n,T_{n+1}],\W11)} 
    \norma{r}_{\L\infty([T_n,T_{n+1}];\L1)}\, .
  \end{eqnarray*}
  The next question is to wonder if $(T_n)$ goes up to $T(h^*)$. We
  assume that it is not the case: then necessarily, $F_{h,n}(T_{n+1})
  \stackrel{n\to \infty}{\longrightarrow} 0,$ since $T_{n+1}-T_n\to
  0$. This is a contradiction to $F_{h,n}(T_{n+1}) = 1/2$.

  Consequently, $T_n\stackrel{n\to \infty}{\longrightarrow}\infty$ and
  the G\^ateaux differentiability is valid for all time $t\in
  [0,T(h^*)]$. Then, making $h^*$ goes to $0$, we obtain that the
  differentiability is valid on in the interval $[0,T(\alpha,\beta)[$.

It remains to check that the G\^ateaux derivative is a bounded linear operator, for $t$ and $\rho_o$ fixed. The linearity is immediate. Additionally, due to the $\L1$ estimate on the solution $r$ of the
  linearized equation (\ref{eq:linear}) given by Proposition
  \ref{prop:existence_r}, we obtain
  \begin{displaymath}
    \norma{ DS_t (\rho_o) (r_o) }_{\L1} 
    = 
    \norma{r(t)}_{\L1}
    \leq  
    e^{Kt\norma{\rho}_{\L\infty(I;\W11)}} e^{Ct}
    \norma{r_o}_{\L1}
    \,,
  \end{displaymath}
  so that the G\^ateaux derivative is bounded, at least for $t\leq
  T<T_{\mathrm{ex}}$.
\end{proofof}

\subsection{Proofs Related to Sections~\ref{sec:Supply}
  and~\ref{sec:Ped}}
\label{subsec:ApplProofs}

\begin{proofof}{Proposition~\ref{prop:OptimalSupply}}
  Note that $v(\rho)$ is constant in $x$, hence $\div V(\rho) = 0$,
  and~\textbf{(A)} is satisfied. Besides, we easily obtain
  $\norma{\pt_x V(\rho)}_{\L\infty(\reali;\reali)} = 0$, $\norma{\pt_x
    V(\rho)}_{\L1(\reali;\reali)}= 0$, $\norma{\pt_x^2
    V(\rho)}_{\L1(\reali;\reali)} = 0$ and
  \begin{eqnarray*}
    \norma{V(\rho_1) - V(\rho_2)}_{\L\infty(\reali;\reali)}
    & \leq &
    \norma{v'}_{\L\infty(\reali;\reali)} \, 
    \norma{\rho_1 - \rho_2}_{\L1(\reali;\reali)}\,,
    \\
    \norma{\pt_x V(\rho_1) - \pt_x V(\rho_2)}_{\L1(\reali;\reali)}
    & = &
    0\,,
  \end{eqnarray*}
  so that~\textbf{(V1)} is satisfied. Similarly, $\pt^2_x V(\rho)=0$
  and $\pt^3_x V(\rho)=0$ imply easily that~\textbf{(V2)}
  and~\textbf{(V3)} are satisfied.

  We consider now~\textbf{(V4)}: is $v$ is $\C2$ then, for all $A, B
  \in \reali$,
  \begin{displaymath}
    v(B) 
    = 
    v(A) 
    + 
    v'(A)(B-A) 
    + 
    \int_0^1 v'' \left(sB + (1-s) A \right)(1-s) (B-A)^2 \d{s}\, .
  \end{displaymath}
  Choosing $A = \int_0^1 \rho(\xi) \d{\xi}$ and $B = \int_0^1 \tilde
  \rho(\xi) \d{\xi}$, we get
  \begin{eqnarray*}
    & &
    \norma{
      v\left(\int_0^1 \tilde\rho(\xi) \d{\xi} \right)
      -
      v\left(\int_0^1 \rho(\xi) \d{\xi} \right)
      -
      v'\left(\int_0^1 \rho(\xi) \d{\xi} \right)
      \int_0^1  \left(\tilde\rho-\rho\right) (\xi)  \d{\xi} 
    }_{\L\infty}
    \\
    & \leq &
    \frac{1}{2}\norma{v''}_{\L\infty}\norma{\tilde\rho-\rho}_{\L1}^2
  \end{eqnarray*}
  and we choose $K = \frac{1}{2}\norma{v''}_{\L\infty}$, $DV(\rho)(r)
  = v'\left(\int_0^1 \rho(\xi) \d{\xi} \right) \int_0^1 r(\xi)
  \d{\xi}$.  Condition~\textbf{(V4)} is then satisfied since there is
  no $x$-dependance, so
  \begin{eqnarray*}
    \norma{V(\tilde\rho)-V(\rho)-DV(\rho)(\tilde\rho-\rho)}_{\W2\infty}
    & = &
    \norma{V(\tilde\rho)-V(\rho)-DV(\rho)(\tilde\rho-\rho)}_{\L\infty}
    \\    
    & \leq & 
    \frac{1}{2}\norma{v''}_{\L\infty}\norma{\tilde\rho-\rho}_{\L1}^2\, .
  \end{eqnarray*}
  Similarly, $\norma{DV(\rho)(r)}_{\W2\infty} =
  \norma{DV(\rho)(r)}_{\L\infty} \leq \norma{v'}_{\L\infty}
  \norma{r}_{\L1}$.  Finally, consider~\textbf{(V5)}:
  \begin{eqnarray*}
    \norma{\div\left[
        v \left( \int_0^1 \tilde\rho(\xi) \d{\xi} \right)
        -
        v \left(\int_0^1 \rho(\xi) \d{\xi} \right)
        -
        v' \left(\int_0^1 \rho(\xi) \d{\xi} \right)
        \int_0^1  \left(\tilde\rho-\rho\right) (\xi) \d{\xi} \right]
    }_{\L1}
    =
    0 \,,
    \\
    \norma{
      \div \left(v'\left(\int_0^1 \rho(\xi)\mathrm{d}\xi\right)
        \int_0^1r(\xi)\mathrm{d}\xi \right)
    }_{\L1}
    =
    0 \,.
  \end{eqnarray*}
  Concluding the proof.
\end{proofof}

\begin{proofof}{Proposition~\ref{prop:ped}}
  The proof exploits the standard properties of the convolution.

  Consider first~\textbf{(V1)}:
  \begin{eqnarray*}
    \norma{\nabla_x V(\rho)}_{\L\infty}
    & = &
    \norma{v'}_{\L\infty} \, \norma{\rho}_{\L\infty}
    \norma{\nabla_x\eta}_{\L1} \,  \norma{\vec{v}}_{\L\infty}
    +
    \norma{v}_{\L\infty} \, \norma{\nabla_x\vec{v}}_{\L\infty}
    \\
    & \leq &
    C(\norma{\rho}_{\L\infty})\, ,
    \\
    \norma{\nabla_x
      V(\rho)}_{\L1}
    & \leq &
    \norma{v}_{\W1\infty}
    \norma{\vec{v}}_{\W11}(1+\norma{\rho}_{\L\infty}
    \norma{\nabla_x\eta}_{\L1})\,,
    \\
    \norma{\nabla_x^2 V(\rho)}_{\L1}
    & \leq &
    \norma{v}_{\W2\infty} \, \norma{\vec{v}}_{\W21}
    \\
    & &
    \times
    \left[
      1 
      + 
      \norma{\rho}_{\L\infty}^2
      \norma{\nabla_x\eta}_{\L1}^2 
      + 
      \norma{\rho}_{\L\infty} 
      \norma{\nabla_x^2\eta}_{\L1} 
      + 
      2 \norma{\rho}_{\L\infty}
      \norma{\nabla_x\eta}_{\L1}
    \right]
    \\
    & \leq &
    C (\norma{\rho}_{\L\infty})\, ,
    \\
    \norma{V(\rho_1) - V(\rho_2)}_{\L\infty}
    & \leq &
    \norma{v'}_{\L\infty} \, \norma{\vec{v}}_{\L\infty} \,
    \norma{\eta}_{\L\infty} \, \norma{\rho_1 - \rho_2}_{\L1}\, ,
    \\
    \norma{\nabla_x \left(V(\rho_1) - V(\rho_2)\right)}_{\L1}
    & = &
    \norma{v}_{\W2\infty}
    \norma{\vec{v}}_{\W1\infty}
    \norma{\eta}_{\W11}
    \left(
      2 + \norma{\nabla_x \eta}_{\L1} \, \norma{\rho_1}_{\L\infty}
    \right)
    \norma{\rho_1-\rho_2}_{\L1} .
  \end{eqnarray*}
  Then, we check~\textbf{(V2)}:
  \begin{eqnarray*}
    \norma{\nabla_x^2 V(\rho)}_{\L\infty}
    & \leq &
    2\norma{v}_{\W2\infty} \, \norma{\vec{v}}_{\W2\infty}
    \\
    & &
    \times
    \left( 
      1 
      + 
      \norma{\rho}_{\L\infty}^2
      \norma{\nabla_x\eta}_{\L1}^2 
      + 
      \norma{\rho}_{\L\infty} 
      \norma{\nabla_x^2\eta}_{\L1} 
      + 
      \norma{\rho}_{\L\infty}
      \norma{\nabla_x\eta}_{\L1}
    \right)\, .
  \end{eqnarray*}
  Entirely analogous computations allow to prove also~\textbf{(V3)}.

  Consider~\textbf{(V4)}. First we look at the Fr\'echet derivative of
  $V(\rho)$: $v$ being $\C2$, we can write, for all $A,B\in \reali$,
  \begin{displaymath}
    v(B)=v(A)+v'(A)(B-A)+\int_0^1 v''(sB+(1-s)A)(1-s)(B-A)^2\,
    \d{s}\, .
  \end{displaymath}
  If we take $A = \rho \ast \eta$ and $B = \tilde\rho \ast \eta$, then
  we get, for $\rho,\tilde \rho \in \L1(\reali^N;\reali)$
  \begin{eqnarray*}
    \norma{\left(
        v \left(\tilde\rho\ast\eta\right)
        -
        v \left(\rho\ast\eta\right)
        -
        v'\left(\rho\ast\eta\right)
        \left((\tilde\rho-\rho)\ast\eta\right)\right)\vec{v}
    }_{\L\infty}
    &\leq &
    \frac{1}{2}\norma{v''}_{\L\infty}
    \norma{\eta}_{\L\infty}^2
    \norma{\tilde\rho-\rho}_{\L1}^2\norma{\vec{v}}_{\L\infty}\,;
  \end{eqnarray*}
  and
  \begin{eqnarray*}
    & &  
    \norma{\nabla_x
      \left[ \left(
          v \left(\tilde\rho\ast\eta\right)
          -
          v \left(\rho\ast\eta\right)
          -
          v'\left(\rho\ast\eta\right)
          \left((\tilde\rho-\rho)\ast\eta\right)\right)\vec{v}\right]
    }_{\L\infty}
    \\
    &\leq &
    \frac{3}{2}
    \norma{v''}_{\L\infty}
    \norma{\eta}_{\W1\infty}^2
    \norma{\tilde\rho - \rho}_{\L1}^2 
    \norma{\vec{v}}_{\W1\infty} 
    \\
    & &
    +
    \frac{1}{2} \norma{v'''}_{\L\infty} 
    \norma{\eta}_{\L\infty}^2 
    \norma{\rho}_{\L\infty} 
    \norma{\nabla_x\eta}_{\L1}
    \norma{\tilde\rho -\rho}_{\L1}^2
    \norma{\vec{v}}_{\L\infty} \,;
    \\
    & & \norma{ \nabla_x^2\left[\left( v
          \left(\tilde\rho\ast\eta\right) - v
          \left(\rho\ast\eta\right) - v'\left(\rho\ast\eta\right)
          \left((\tilde\rho-\rho)\ast\eta\right) \right)\vec{v}\right]
    }_{\L\infty}
    \\
    & \leq & \norma{v^{(4)}}_{\L\infty}
    \norma{\tilde\rho-\rho}_{\L1}^2 \norma{\eta}_{\L\infty}^2
    \norma{\nabla_x\eta}_{\L1}^2 \left( \norma{\rho}_{\L\infty} +
      \norma{\tilde \rho}_{\L\infty} \right)^2
    \norma{\vec{v}}_{\L\infty}
    \\
    & & + 2 \norma{v^{(3)}}_{\L\infty} \norma{\tilde\rho-\rho}_{\L1}^2
    \norma{\eta}_{\W1\infty}^2 \left( \norma{\rho}_{\L\infty} +
      \norma{\tilde\rho}_{\L\infty} \right) \norma{\nabla \eta}_{\L1}
    \norma{\vec{v}}_{\L\infty}
    \\
    & & + \frac{1}{2} \norma{v^{(3)}}_{\L\infty}
    \norma{\tilde\rho-\rho}_{\L1}^2 \norma{\eta}_{\L\infty}^2 \left(
      \norma{\rho}_{\L\infty} \norma{\eta}_{\W21} +1 \right)
    \norma{\vec{v}}_{\W1\infty}
    \\
    & & + 6 \norma{v''}_{\L\infty} \norma{\tilde\rho-\rho}_{\L1}^2
    \norma{\eta}_{\W2\infty}^2 \norma{\vec{v}}_{\W2\infty} \,.
  \end{eqnarray*}
  Then, $DV(\rho)(r) = v' ( \rho \ast \eta) r\ast\eta \, \vec{v}$.

  In order to satisfy~\textbf{(V4)}, we have also to check that the
  derivative is a bounded operator from $\C2$ to $\L1$. We have,
  \begin{eqnarray*}
    \norma{ DV(\rho)(r)}_{\L\infty} 
    & \leq & 
    \norma{v'}_{\L\infty}
    \norma{\eta}_{\L\infty}
    \norma{\vec{v}}_{\L\infty} 
    \norma{r}_{\L1}\,, 
    \\
    \norma{\nabla_x DV(\rho)(r)}_{\L\infty} 
    & \leq & 
    \norma{v}_{\W2\infty} 
    \norma{\vec{v}}_{\W1\infty} 
    \norma{\eta}_{\W1\infty}
    \left(
      2
      +
      \norma{\rho}_{\L\infty } 
      \norma{\eta}_{\W11}
    \right)
    \norma{r}_{\L1} \,, 
    \\
    \norma{\nabla^2_x DV(\rho)(r)}_{\L\infty} 
    & \leq & 
    \norma{v}_{\W3\infty} 
    \norma{\eta}_{\W2\infty}
    \norma{\vec{v}}_{\W2\infty}
    \\
    & &
    \times
    \left(
      4
      +
      5 
      \norma{\rho}_{\L\infty}
      \norma{\eta}_{\W21} 
      +
      \norma{\rho}_{\L\infty}^2
      \norma{\nabla\eta}_{\L1}^2 
    \right)
    \norma{r}_{\L1}\,.
  \end{eqnarray*}
  Finally, we check that also~\textbf{(V5)} is satisfied:
  \begin{eqnarray*}
    & &
    \norma{\div 
      \left(
        V(\tilde\rho) - V(\rho) - DV(\rho)(\tilde\rho-\rho) 
      \right)
    }_{\L1}
    \\
    & \leq &  
    \frac{1}{2} \norma{v}_{\W3\infty} 
    \norma{\tilde\rho-\rho}_{\L1}^2 
    \norma{\eta}_{\L1} \norma{\eta}_{\W1\infty} 
    \norma{\vec{v}}_{\W1\infty} 
    \left(
      3
      + 
      \norma{\rho}_{\L\infty}
      \norma{\eta}_{\W11} 
    \right)\, ,
    \\
    & & \norma{\div DV(\rho)(r)}_{\L1}
    \\
    & = & \norma{\div (v'(\rho\ast\eta)r\ast \eta\, \vec{v})}_{\L1}
    \\
    & \leq & \norma{v}_{\W2\infty} \norma{\eta}_{\W11}
    \norma{\vec{v}}_{\W1\infty} \left(2+\norma{\rho}_{\L\infty}
      \norma{\nabla_x\eta}_{\L1} \right) \norma{r}_{\L1}
  \end{eqnarray*}
  completing the proof.
\end{proofof}

\begin{remark}
  The above proof shows that condition~\textbf{(B)} is not satisfied
  by~(\ref{eq:Ped}). Indeed, here we have that $C$ grows linearly:
  $C(\alpha) = 1+\alpha$. Hence, with the notation used in the proof
  of Theorem~\ref{thm:V0}, for $\alpha_1 >0$, we have
  \begin{displaymath}
    \sum_{k=1}^n T( \alpha_k,\alpha_{k+1})
    \leq 
    \sum_{k=1}^n \frac{1}{1+\alpha_{k+1}} 
    \int_{\alpha_k}^{\alpha_{k+1}} \frac{1}{t} \d{t}
    \leq 
    \sum_{k=1}^n 
    \int_{\alpha_k}^{\alpha_{k+1}} \frac{1}{(1+t)t} \d{t}
    \leq
    \int_{\alpha_1}^{+\infty} \frac{1}{(1+t)t} \d{t}
  \end{displaymath}
  and the latter expression is bounded. This shows that, in the case
  of~(\ref{eq:Ped}), the technique used in Theorem~\ref{thm:V0} does
  not apply.
\end{remark}
\medskip

\noindent\textbf{Acknowledgment.}
{\small{  This work has been supported by the DFG SPP1253,  
DAAD  D/06/19582, DAAD D/08/11076, and HE5386/6-1. The third author
 thanks Sylvie Benzoni-Gavage for useful discussions.
}}

\bibliographystyle{abbrv}

\bibliography{controlN}

\end{document}